\documentclass[times]{MyArticle}

\usepackage[colorlinks,bookmarksopen,bookmarksnumbered,citecolor=red,urlcolor=red]{hyperref}

\usepackage[usenames,dvipsnames,svgnames,table]{xcolor}

\newcommand\BibTeX{{\rmfamily B\kern-.05em \textsc{i\kern-.025em b}\kern-.08em
T\kern-.1667em\lower.7ex\hbox{E}\kern-.125emX}}

\begin{document}

\title{A Meshfree Generalized Finite Difference Method for Surface PDEs}

\runningheads{Suchde \textit{et~al.}}{A Meshfree GFDM for Surface PDEs}

\author{Pratik Suchde \affil{1}\corrauth, J\"org Kuhnert\affil{1}}

\address{\affilnum{1}Fraunhofer ITWM, 67663 Kaiserslautern, Germany}

\corraddr{P. Suchde. E-mail: pratik.suchde@itwm.fraunhofer.de}

\begin{abstract}
In this paper, we propose a novel meshfree Generalized Finite Difference Method~(GFDM) approach to discretize PDEs defined on manifolds. Derivative approximations for the same are done directly on the tangent space, in a manner that mimics the procedure followed in volume-based meshfree GFDMs. As a result, the proposed method not only does not require a mesh, it also does not require an explicit reconstruction of the manifold. In contrast to some existing methods, it avoids the complexities of dealing with a manifold metric, while also avoiding the need to solve a PDE in the embedding space. A major advantage of this method is that all developments in usual volume-based numerical methods can be directly ported over to surfaces using this framework. We propose discretizations of the surface gradient operator, the surface Laplacian and surface Diffusion operators. Possibilities to deal with anisotropic and discontinous surface properties (with large jumps) are also introduced, and a few practical applications are presented. 
\end{abstract}


\maketitle

\vspace{-6pt}

%%%%%%%%%%%%%%%%%%%%%%%%%%%%%%%%%%%%%%%%%%%%%%%%%%%%%%%%%%%%%%%%%%%%%%%%%%%%%%%%%%%%%%
%%%%%%%%%%%%%%%%%%%%%%%%%%%%%%%%%%%%%%%%%%%%%%%%%%%%%%%%%%%%%%%%%%%%%%%%%%%%%%%%%%%%%%
%%%%%%%%%%%%%%%%%%%%%%%%%%%%%%%%%%%%%%%%%%%%%%%%%%%%%%%%%%%%%%%%%%%%%%%%%%%%%%%%%%%%%%
%%%%%%%%%%%%%%%%%%%%%%%%%%%%%%%%%%%%%%%%%%%%%%%%%%%%%%%%%%%%%%%%%%%%%%%%%%%%%%%%%%%%%%

\section{Introduction}

The solution to partial differential equations~(PDEs) defined on a surface or manifold is of fundamental interest in various fields. They have application in the fields of computer graphics~\cite{Turk1991}, image processing~\cite{Diewald2000}, fluid flow~\cite{Myers2002}, and cell biology~\cite{Novak2007}, to name a few.

%A significant amount of work has been done to solve PDEs on surfaces discretized by a mesh using both finite element methods \cite{Dziuk2013} and finite volume methods \cite{Du2006}.  . However, both a parametrization and a meshing of the surface can be non-trivial to obtain. Thus, there is a need to solve surface PDEs directly on a scattered point cloud in a meshfree way. 

Most existing surface PDE solvers can be classified into two types. The first, referred to as intrinsic methods, solve the PDE directly on the manifold using either a mesh~\cite{Du2003,Du2006,Dziuk2013,Olshanskii2009}, a parametrization of the manifold~\cite{Floater2005}, or an explicit reconstruction of the manifold \cite{Liang2013}. In the second type of methods, referred to as embedding methods, the surface PDE is extended to a PDE defined on a band around the manifold, which is then discretized \cite{Bertalmio2001,Chu2018, Ratz2006,Ruuth2008}. The significant disadvantage of these methods is that they rely on discretizing a higher dimensional space, and thus they tend to be expensive in terms of computational time. On the other hand, the intrinsic methods have the disadvantage that the parametrization can be non-trivial to obtain, and that they need to deal with singularities arising in metric terms of surface differential operators~\cite{Fuselier2013}. A more detailed breakdown of existing literature on surface PDEs can be found in \cite{ChenThesis_CPM,IngridThesis_CPM}. 

Constructing a good meshing of a manifold can be a very difficult process. It is thus often desirable to discretize a manifold directly with a set of scattered points, referred to as a point cloud. As a result, the need for meshfree methods to solve surface PDEs arises. In this paper, we present a new meshfree method for solving surface PDEs.

Several meshfree methods for surface PDEs have already been proposed. Most notable among these are the radial-basis function~(RBF) based methods \cite{Flyer2014,Fuselier2013}. RBF based surface PDE solvers combine the advantages of both intrinsic and embedding methods. They solve surface PDEs on the surface itself, without needing to parametrize or mesh the surface. However, they suffer from the drawback of needing a somewhat ad-hoc choice of basis function and related shape parameter, which affects stability and the conditioning of the linear systems \cite{Fuselier2013,Mongillo2011}. This choice is application dependent \cite{Simonenko2014}, given data~(right hand side or initial condition) dependent \cite{Davydov2011} , and can even be domain dependent. The optimal choice is not always known. We use a meshfree generalized finite difference method~(GFDM) to avoid this issue.

Meshfree GFDMs \cite{Fan2018,Gavete2017,Katz2010,Luo2016} are strong form meshfree methods that have been shown to be robust methods, and have been used in a wide variety of applications \cite{Drumm2008,Jefferies2015,Moller2007,Tramecon2013}. Approximations are carried out using a weighted least squares approach. In this paper, we propose a meshfree GFDM approach to solve PDEs defined on a surface. This is done by projecting local neighbourhoods to the tangent space, and performing approximations there. This extends the work of Demanet~\cite{Demanet2006} for mesh-based surface Laplacians to a meshfree context applicable to various differential operators. The proposed new method retains the advantages of RBF based methods, while avoiding the disadvantage of the ad-hoc choice of basis functions and shape parameters.

We note that meshfree GFDM approaches to solve surface PDEs have already been proposed by Liang \textit{et~al.} \cite{Liang2012,Liang2013}. In contrast, we do not rely on surface-based metrics, and thus avoid the issue of arising singularities. Another key difference from \cite{Liang2012,Liang2013} is the formulation used here enables transferring developments from volume-based numerical methods directly to surfaces. 

The remainder of the paper is organized as follows. In Section~\ref{sec:Definitions}, we present the basics of meshfree GFDMs and introduce the notation used in the paper. Section~\ref{sec:SurfaceDOs} introduces a novel way to discrete surface differential operators in a meshfree GFDM setting. Section~\ref{sec:BoundaryConditions} contains a short note on the implementation of boundary conditions. While Section~\ref{sec:SurfaceDOs} deals with $2$-manifolds in $3$ dimensional space,  Section~\ref{sec:HighD_coD} extends those ideas to higher dimensions and co-dimensions. Section~\ref{sec:NumResults} presents a range of numerical examples and validation, and the paper is concluded with a discussion on the work in Section~\ref{sec:Conclusion}.

%
%\section{Basic Definitions}
\section{Preliminaries}
\label{sec:Definitions}

\subsection{Basics and Nomenclature}

For the majority of this paper, we consider a smooth orientable $2$-manifold $M$ embedded in $\mathbb{R}^3$. Extensions of the ideas presented here to higher dimensions and co-dimensions can be done easily, and will be discussed briefly in Section \ref{sec:HighD_coD}. The manifold $M$ is also referred to as the embedded space or surface, while $\mathbb{R}^3$ is also referred to as the embedding space.

We consider a manifold given by an unevenly distributed point cloud consisting of $N$ points.  For each point $i$ on the manifold, approximations are carried out on the neighbourhood or support $S_i$ consisting of $|S_i|$ nearby points. $S_i$ is based on proximity, given by Eucledian distances in the embedding space. Throughout this paper, all distances are computed in the embedding space only. Distances along the manifold are never needed, and are not computed. The size of the support $S_i$ is given by the smoothing length $h_i$. We adopt the following distance conventions from volumetric meshfree GFDM for fluid flow \cite{Drumm2008, Suchde2017_CCC}. During set up of the point cloud, it is ensured that no two points are closer than $r_{min}h$, and that there is no hole of size $r_{max}h$ on the surface that does not contain any point. The parameters $r_{min}=0.2$ and $r_{max}=0.45$ are fixed across all simulations, and are adopted from volumetric meshfree GFDM conventions \cite{Suchde2017_CCC}. This ensures that that each support of size $h$ has about $15-20$ points. This also results in $h$ serving as an indication of the point cloud spacing. In all simulations, it is assumed that $h$ is chosen such that the local neighbourhoods are sufficiently dense to accurately capture the local curvature of the surface. 

We assume that at each interior point $i$, the unit surface normal $\vec{n}_i$ and unit tangents $\vec{t}_{1,i}$, $\vec{t}_{2,i}$ to the surface have already been computed. Similarly, for each boundary point, the unit surface normal $\vec{n}_i$, the unit surface tangent $\vec{t}_i$, and the unit boundary normal $\vec{\nu}_i$ are assumed to be known. Note that these normals and tangents form an orthogonal system of vectors. Normal and tangent computation could have been done in any of multiple ways. One possibility is to use Principle Component Analysis~(PCA) to construct the normal and tangent information based on the eigenvalues of local covariance matrices (see, for example, \cite{Liang2013,Mitra2003}). Weighted PCA approaches \cite{Petronetto2013} could also be used for the same. Alternatively, surface normals and tangents could be computed based on local geometric information in each neighbhourhood using the same procedures done in meshfree GFDM simulations of fluids to compute normals at the free surface (see, for example, \cite{Edgar2017}). This latter approach will be used for all simulations in this paper. Further, if the point cloud is based directly on a CAD model, normal information could also be taken from the CAD model itself, if it is available, and the tangents could be computed accordingly.

The following notation is used for differential operators throughout this paper. $\nabla$, $\Delta$ are used to denote the continuous operators in the embedding space $\mathbb{R}^3$. Subscripts $M$ and $T$ to the differential operators are used to indicate the corresponding operator defined on the manifold and the tangent space respectively, while the subscript $tn$ is used to indicate operators in $\mathbb{R}^3$ in a rotated tangent-normal frame. A $\tilde{}\,$ is used overhead to explicitly state, if needed, that the operator being considered is a discrete one. A second subscript $i$ is used to indicate that the discrete operator is being considered at the point $i$. Thus, $\widetilde{\Delta}_{M,i}$ indicates the discrete surface Laplacian at point $i$.

\subsection{Volume-based Meshfree GFDMs}
\label{sec:VolumetricGFDMs}

Classical meshfree GFDMs for volume or bulk domains shall henceforth be referred to as `volumetric' GFDMs. They are strong-form methods, in which for a function $u$ defined at each numerical point $i=1,2,\dots,N$, its derivatives are approximated as
\begin{align}
	\nabla u(\vec x_i) &\approx \widetilde{\nabla}_i u = 
	\left( \def\arraystretch{2.0} \begin{array}{c}
		\sum_{j\in S_i}c_{ij}^x u_j	\\	
		\sum_{j\in S_i}c_{ij}^y u_j	\\
		\sum_{j\in S_i}c_{ij}^z u_j	\\				
	\end{array}  \right) \,,\label{Eq:GFDM_Definition1}\\
	\Delta u(\vec x_i) &\approx \widetilde{\Delta}_i u = \sum_{j\in S_i}c_{ij}^\Delta u_j\,,\label{Eq:GFDM_Definition2}
\end{align}
%
%
%\begin{equation}
%	\label{Eq:GFDM_Definition}
%	 %\partial^* u(\vec x_i)\approx \tilde{\partial}^*_i u = \sum_{j\in S_i}c_{ij}^*u_j\,,
%
%\end{equation}
%
where $\widetilde{(\cdot)}_i$ indicates the discrete differential operator at point $i$, and the differential operators without the tilde indicate the continous differential operators. All meshfree GFDM notation used here follows from \cite{Suchde2017_CCC}. For each point $i$, the stencil coefficients $c_{ij}$ are found using a weighted least squares approach. The weighted sum of the stencil coefficients is minimized such that monomials $m\in\mathcal{P}$ up to a certain order, usually $2$, are exactly differentiated. For example, for the Laplacian, 
\begin{align}
	\sum_{j\in S_i}c_{ij}^{\Delta}m_j &= \Delta m (\vec{x_i}) \qquad \forall m\in\mathcal{P}\,,\label{Eq:Consistency}\\
	\text{min } J_i &= \sum_{j\in S_i} \left( \frac{c_{ij}^{\Delta}}{W_{ij}} \right)^2\,, \label{Eq:BasicMin}
\end{align}
where $W$ is a weighting function. Throughout this paper, we use a Gaussian weighting function
\begin{equation}
	W_{ij}=\exp(-W_F \frac{\|\vec{x}_j-\vec{x}_i\|^2}{h_i^2 + h_j^2})\,,
\end{equation}
where $W_F>0$, and the inclusion of the smoothing length $h_j$ of the neighbour point is for point clouds with uneven point distributions. Note that for a central point $i$, the weighting function is only relevant for the neighbouring points $j\in S_i$. The same procedure also holds for other differential operators. The monomials depend on distance relative to the central point $i$. So, for example, in $\mathbb{R}^2$, for the $x$ derivative, the monomial consistency conditions lead to the following system
\begin{equation}
\label{Eq:DO_Example}
\underbrace{\left(\begin{array}{ccc}
	\cdots& 1 &\cdots\\
	\cdots& \delta x_{ij} &\cdots\\
	\cdots& \delta y_{ij} &\cdots\\
	\cdots& \delta x_{ij}^2 &\cdots\\
	\cdots& \delta y_{ij}^2 &\cdots\\
	\cdots& \delta x_{ij}\delta y_{ij} &\cdots\\
	\end{array}\right)}_{M_i^T}
\underbrace{\left(\begin{array}{c}	
	\vdots \\
	c_{ij}^x\\
	\vdots\\
	\end{array}\right)}_{\vec{c}_i^{\,x}} = 
\underbrace{\left(\begin{array}{c}
	0\\
	1\\
	0\\
	0\\
	0\\
	0\\
	\end{array}\right)}_{\vec{b}_i^{\,x}}\,,
\end{equation}
where $\delta x_{ij}=x_j - x_i$ and $\delta y_{ij} = y_j - y_i$. The number of neighbours is always taken to be larger than the number of monomial consitency conditions. Continuing the example of the $x$ derivative, the minimization equivalent of Eq.\,\eqref{Eq:BasicMin} can be written as $\text{min } J_i = \|W_i^{-1}\vec{c}_i^{\,x}\|^2$, with $W_i$ being the diagonal matrix of weights. This minimization leads to the differential operators given by $\vec{c}_i^{\,x} = W_i^2M_i\left(M_i^TW_i^2M_i\right)^{-1}\vec{b}_i^{\,x}$. Note that the LHS matrix in Eq.\,\eqref{Eq:DO_Example} is the same for all differential operators being approximated. Thus, multiple differential operators can be computed numerically with a single minimization procedure with multiple right hand sides~\cite[Section A.2]{Suchde2018_Thesis}. %Numerically, rather than computing the numerical inverse, the minimization problem is solved with a QR decomposition, 

We wish to use a similar approach here for surface derivatives. All numerical differential operators on a manifold will be defined in a manner similar to Eq.\,\eqref{Eq:GFDM_Definition1} and Eq.\,\eqref{Eq:GFDM_Definition2}, and the procedure for computing the stencil coefficients will be done in a similar manner to Eq.\,\eqref{Eq:Consistency} and Eq.\,\eqref{Eq:BasicMin}.

\section{Surface Differential Operators using a Meshfree GFDM}
\label{sec:SurfaceDOs}

Consider a function $u: M \rightarrow \mathbb{R}$ defined on the surface, and an open subset $\Omega$ of $\mathbb{R}^3$ containing $M$ i.e. $M \subset \Omega \subset \mathbb{R}^3$. A function $\hat{u}: \Omega \rightarrow \mathbb{R}$ is said to be an extension of $u$ if $u$ and $\hat{u}$ agree on the manifold, $\hat{u}|_{M} = u$. Further, $\hat{u}$ is said to be a normal extension of $u$ if $\vec{n}\cdot\nabla \hat{u} \equiv 0$.

In this paper, we use normal extensions to extend functions \textit{locally}. We then use these extensions to define numerical differential operators entirely on the tangential plane at each point. We emphasize that the extension is \textit{not} done numerically as is the case for embedding methods such as the closest point method~\cite{Marz2012,Ruuth2008}. Further, no dummy or virtual points are needed for differential operator compuation, as is the case for some RBF-based methods for surface PDEs~\cite{Piret2012}. Here, only the manifold $M$ is discretized, whereas $\Omega$ is not. The concept of (normal) extensions is only used to derive a way to compute differential operators directly on the tangential plane in a straight forward manner. 

\subsection{Surface Gradient}
\label{sec:SurfaceGradient}

The surface gradient of a function can be defined as the conventional gradient of its extension with the component normal to the surface removed
\begin{align}
	\nabla_{M} u &= \nabla \hat{u} - \vec{n} \left( \vec{n} \cdot \nabla \hat{u} \right)\,, \label{Eq:SurfaceGradient_Basic}\\
	&= \left( \mathbf{P}\nabla \right) \hat{u}\,,
\end{align}
where $\mathbf{P}=\mathbf{I} - \vec{n}\vec{n}^T$ for identity $\mathbf{I}$ is the projection matrix (to the tangent space). Here, we define $\hat{u}$ such that it is a normal extension of $u$. Thus, $\vec{n}\cdot\nabla \hat{u}= 0$, and Eq.\,\eqref{Eq:SurfaceGradient_Basic} reduces to
\begin{equation}
	\label{Eq:GradientTheorem}
	\nabla_{M} u = \nabla \hat{u} \,.
\end{equation}
Thus, determining a numerical approximation $\widetilde{\nabla}\hat{u}$ to $\nabla\hat{u}$  (the gradient of the extension $\hat{u}$ ),  will give us an approximation to the surface gradient $\widetilde{\nabla}_{M}u$ of the original function $u$. 

To compute the approximation $\widetilde{\nabla}\hat{u}$, we once again make use of the fact that $\vec{n}\cdot\nabla \hat{u}= 0$ everywhere on the manifold. Rather than computing $\widetilde{\nabla}\hat{u}$ directly, we first consider the rotated components
\begin{equation}
	\nabla_{tn}\hat{u} = \left(\begin{array}{c} 
	\vec{t}_1\cdot\nabla \hat{u} \\
	\vec{t}_2\cdot\nabla \hat{u} \\
	\vec{n}\cdot\nabla \hat{u} \\	
	\end{array}\right) = \left(\begin{array}{c} 
	\vec{t}_1\cdot\nabla \hat{u} \\
	\vec{t}_2\cdot\nabla \hat{u} \\
	0 \\	
	\end{array}\right) \,,
\end{equation}
where the subscript $tn$ indicates the rotated $\vec{t}_1$, $\vec{t}_2$, $\vec{n}$ coordinate frame at a particular point. Thus, we only need to compute numerical approximations to the tangential components $\vec{t}_k \cdot \nabla \hat{u}$ for $k=1,2$. Once these tangential components are known, they can be rotated to get the gradient $\nabla \hat{u}$.
\begin{equation}
	\label{Eq:RotationOfGradient}
	\nabla \hat{u} =  \underbrace{ \left(\begin{array}{ccc} 
	\vec{t}_1 & \vec{t}_2 & \vec{n} 
	\end{array}\right)}_{R^T}      
	\nabla_{tn}\hat{u} \,,
\end{equation}
where the rotation matrix $R^T$ is composed of the the tangential and normal column vectors. 

\subsubsection{Numerical Surface Gradient Operator}
~\\
We showed above that the problem of computing numerical differential operators for the surface gradient of a function can be reduced to that of computing numerical differential operators for the tangential components of the regular volumetric gradient of the normal extension of the function.

Numerically, at each discrete point $i$ on the manifold, we compute approximations of the tangential components of $\nabla \hat{u}$ entirely on the tangent plane $T_i$ spanned by $\vec{t}_{1,i}$ and $\vec{t}_{2,i}$, i.e. we compute
\begin{equation}
	\widetilde{\nabla}_{T} u = \left(\begin{array}{c} 
	\vec{t}_1\cdot \widetilde{\nabla} \hat{u} \\
	\vec{t}_2\cdot \widetilde{\nabla} \hat{u} 
	\end{array}\right) \,,
\end{equation}
where $\widetilde{\nabla}_{T}$ represents the discrete $2$ dimensional gradient in the tangential plane. For this, we project each point $j \in S_i$ to the tangential plane $T_i$. Let the projection of the point $j \in S_i$ to $T_i$ be labelled as $j_{T_i}$. Since $u$ is being extended normally (the derivative in the normal direction is $0$), $\hat{u}$ evaluated at point $j_{T_i}$  is the same as $u$ evaluated at point $j$. Thus, we set
\begin{equation}
	\label{Eq:Projection_j_Equiv}
	\hat{u}_{j_{T_i}} = u_j
\end{equation}
A neighbouring point $j \in S_i$ is projected to the tangent plane $T_i$ along the surface normal of the central point $\vec{n}_i$. As a result, tangential distances are preserved. The distances between the central point $i$ and the projected points $j_{T_i}$ can be computed by simply rotating the original distances. If $\delta\vec{x}_{ij} = \vec{x}_j - \vec{x}_i$, then the distances in the tangential and normal directions are given simply by $R\, \delta\vec{x}_{ij}$, where $R$ is the rotation matrix introduced in Eq.\,\eqref{Eq:RotationOfGradient}. Of these, only the tangential distances are needed for derivative computation. The computation of the differential operators according to a procedure like Eq.\,\eqref{Eq:Consistency} and Eq.\,\eqref{Eq:BasicMin} only requires the distances between points, as shown in Eq.\,\eqref{Eq:DO_Example}. Since we already have these distances, there is no need to actually add a point numerically at the projected location. 

\begin{figure}
  \centering
  \includegraphics[width=0.45\textwidth]{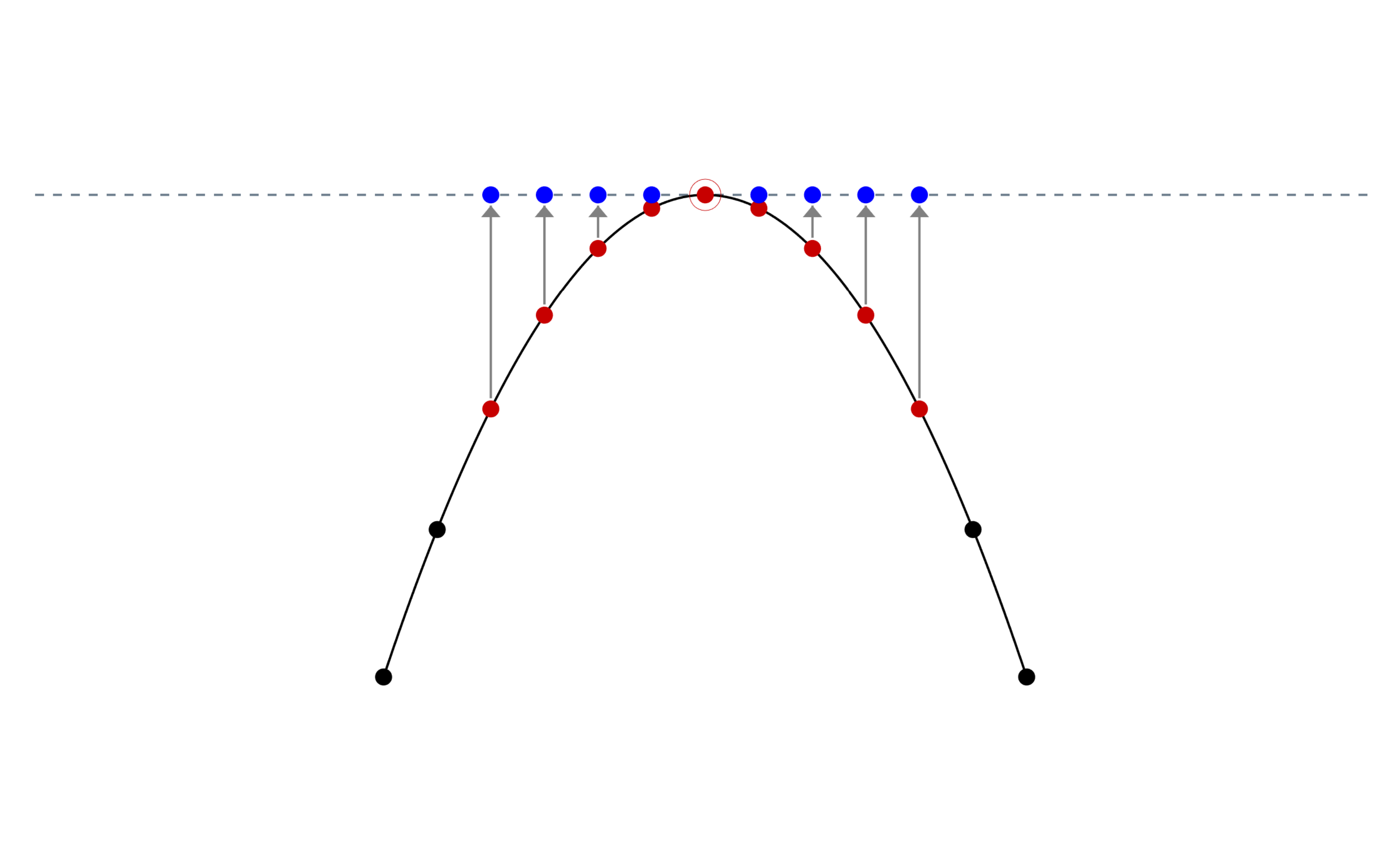}
  \includegraphics[width=0.45\textwidth]{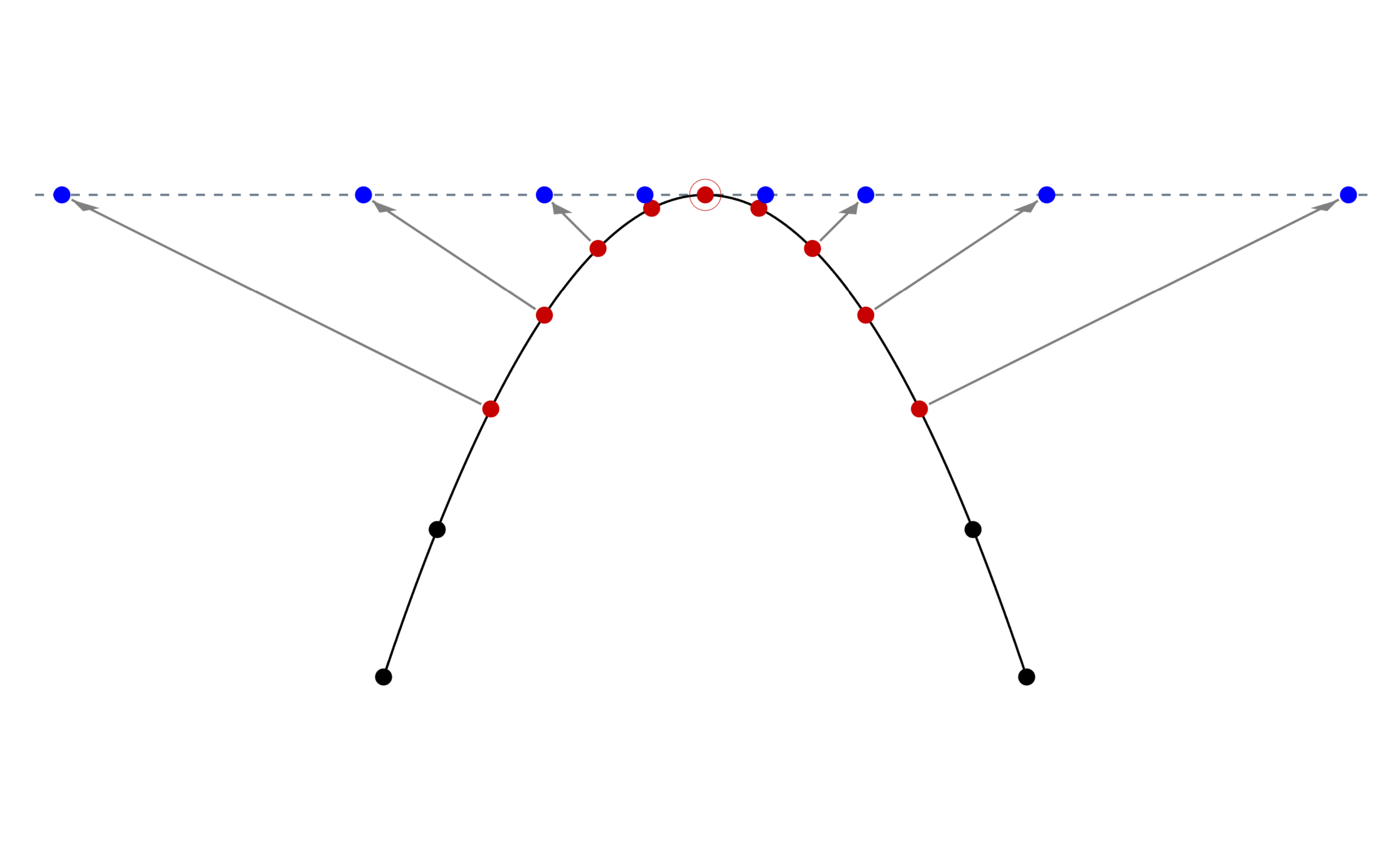}
  \caption{Projection of neighbouring points to the tangent space along the central normal~(left), and the neighbour normals~(right). The central point is shown with an additional circle around it, and all its neighbouring points are marked in red. The remaining black points on the manifold are shown for reference. The projected locations on the tangent space are marked in blue. The manifold is shown in black, while the tangent line at the central point is shown with a blue dashed line.}%
  % Figure created in Paraview - see SaveState3_forDomain.pvsm in Manifold_Results/FiveStripProblem/
  \label{Fig:Projection}
\end{figure}

We note that the idea of projecting to the tangent space to compute numerical approximations is not a novel one. It has also been done by mesh-based surface PDE solvers. Lai \textit{et~al.} \cite{Lai2013} do the same, but while also working with surface-based metrics. Demanet~\cite{Demanet2006} does the same in the mesh-based framework, but only for the surface Laplacian. There, they project points along the surface normal of the neighbouring point itself $\vec{n}_j$. Such a procedure would involve a little more work numerically, and as we show later, is less accurate than the central normal projection used here. The difference between the two methods of projecting neighbouring points to the tangent space is illustrated in Figure~\ref{Fig:Projection} for a $1$-dimensional manifold in $\mathbb{R}^2$.
% Does Lai2013 do central normal projection?

Unlike the neighbour normal projection, the central normal projection method does not give a global normal extension of the function. i.e. if the union of all the projected locations on tangential planes $\mathcal{T} = \cup_{i,j\in S_i}\, \vec{x}_{j_{T_i}}$ is considered, 
$\hat{u}$ defined on $\mathcal{T}$ does not form a global normal extension of $u$. However, this is not relevant in the present context. Our interest is only local. For each point $i$, a virtual projection of its neighbours in $S_i$ to its tangential plane $T_i$ is only performed to numerically compute the derivatives at point $i$ itself. 
%Since function extensions are only being performed locally, and not globally, projection neighbouring point along the central point $i$ normal is not wrong. 
%

For each point $i$, once the tangential distances to the neighbouring points are known, we  compute $2$ dimensional volumetric numerical differential operators on the tangential plane for the first derivatives along the $\vec{t}_{1,i}$ and $\vec{t}_{2,i}$ directions, as explained in Section~\ref{sec:VolumetricGFDMs}.
\begin{align}
	\sum_{j\in S_i}c_{ij_{T}}^{t_k}m_{j_{T}} &= \frac{\partial}{\partial t_k} m (\vec{x}_i) \qquad \forall m\in\mathcal{P}_{T}\,,\label{Eq:TP_Consistency}\\
	\text{min } J_i &= \sum_{j\in S_i} \left( \frac{c_{ij_{T}}^{t_k}}{W_{ij_{T}}} \right)^2\,, \label{Eq:TP_Min}
\end{align}
for $k=1,2$, where $\mathcal{P}_{T}$ are the set of monomials, usually up to order $2$, in $\vec{t}_{1,i}$ and $\vec{t}_{2,i}$ on the tangent plane. 

Using Eq.\,\eqref{Eq:RotationOfGradient} and Eq.\,\eqref{Eq:GradientTheorem}, the computed differential operator stencil coefficients can now be rotated to obtain the numerical surface gradient operator. We have
\begin{equation}
	\label{Eq:SurfGrad_FullDefinition}
	\widetilde\nabla_{M,i} u = %
									\left(\def\arraystretch{2.0} \begin{array}{c}
									\sum_{j \in S_i}c_{ij}^{M,x} u_j \\												\sum_{j \in S_i}c_{ij}^{M,y} u_j \\		
									\sum_{j \in S_i}c_{ij}^{M,z} u_j \\		
										\end{array}\right)\,,
\end{equation}
where $c_{ij}^{M,x}$ are the stencil coefficients for the surface gradient in the $x$ direction, and similarly for the other directions. Further, 
\begin{equation}
	\left(\def\arraystretch{2.0} \begin{array}{c}
	c_{ij}^{M,x} \\														   			c_{ij}^{M,y} \\		
	c_{ij}^{M,z} \\		
		\end{array}\right) = %
	R^T 		\left(\def\arraystretch{2.0} \begin{array}{c}
	c_{ij}^{t_1} \\														   						c_{ij}^{t_2} \\		
	c_{ij}^{n} \\		
		\end{array}\right)\,,
\end{equation}
for $c_{ij}^n = 0$.

This procedure used for the computation of the numerical surface gradient operators can be extended easily to obtain discretizations for any surface differential operator. 

A key point to note here is that the main computation of differential operator stencil coefficients, according to Eq.\,\eqref{Eq:TP_Consistency} and Eq.\,\eqref{Eq:TP_Min}, follows the exact same procedure as that of regular volumetric differential operator computation. As a result, any and all developements in volumetric meshfree GFDMs can be directly used for surface PDEs. This is the one of the biggest advantages of this method over existing work for the same. A significant amount of work has been done to modify GFDM differential operator for different ends, and they can be easily carried over to surface PDEs using the present work. These include developments such as higher order spatial discretizations~\cite{Milewski2012}, conservation~\cite{Chiu2012, Suchde2017_CCC}, accuracy considerations~\cite{Suchde2018_INSE}, upwinding methods for advection \cite{Praveen2007, SeifarthThesis}, staggered methods \cite{Trask2017}, among others. A few examples of carrying over volumetric GFDMs developments to surface PDEs are shown in the coming sections.

%In summary, the surface gradient of the function defined on the manifold is the same as the (regular) gradient of a normal extension of the function. This is evaluated in a rotated tangent-normal frame of reference. The normal component is $0$ by definition, and the tangential components are computed numerically by the gradient on the tangential plane at each numerical point.\\

%$\hat{x}=(1,0,0)$, $\hat{y}=(0,1,0)$ and $\hat{z}=(0,0,1)$ are the standard basis vectors in $\mathbb{R}^3$.

\subsection{Surface Divergence}
\label{sec:SurfaceDivergence}

Consider a vector valued function $\vec{v} = (v^1, v^2, v^3)$ defined on the manifold. A normal extension of $\vec{v}$ is obtained by a normal extension of each component of $\vec{v}$. i.e. $\vec{\hat{v}} = \left( \hat{v}^1, \hat{v}^2, \hat{v}^3 \right)$. Now, the surface divergence of $\vec{v}$ can be written as
\begin{align}
	\nabla_{M}\cdot \vec{v} &=  \left(\mathbf{P}\nabla\right) \cdot \vec{\hat{v}}\,,\\
	&= \nabla\cdot\vec{\hat{v}} - \vec{n} \cdot \left[ \left( \vec{n} \cdot \nabla \right) \vec{\hat{v}} \right] \,,\\
	&= \nabla\cdot\vec{\hat{v}} \label{Eq:DivergenceTheorem}\,.
\end{align}
The numerical gradient operators computed above can be used to compute the surface divergence. For brevity, we write the numerical surface gradient according to Eq.\,\eqref{Eq:SurfGrad_FullDefinition} as
\begin{equation}
	\label{Eq:SurfaceGradientConcise}
	\widetilde{\nabla}_{M,i}u = ( G_{1i} u, G_{2i} u, G_{3i} u)^T\,,
\end{equation}
where $G_{ki}$ are the operators computed in the previous section. For example, $G_{1i}u = \sum_{j\in S_i}c_{ij}^{M,x}u_j$. Using the notation of Eq.\,\eqref{Eq:SurfaceGradientConcise}, the numerical surface divergence can be computed as
\begin{equation}
	\widetilde{\nabla}_{M}\cdot \vec{v} = \sum_{k=1,2,3} G_{ki} v^k\,.
\end{equation}

\subsection{Surface Laplacian}
\label{sec:SurfaceLaplace}

The surface Laplacian, or Laplace Beltrami, of a scalar valued function is defined as
\begin{equation}
	\label{Eq:LapBeltDefinition}
	\Delta_{M} u = \nabla_{M}\cdot \nabla_{M} u\,.
\end{equation}
Using Eq.\,\eqref{Eq:GradientTheorem} and Eq.\,\eqref{Eq:DivergenceTheorem} leads to 
\begin{equation}
	\label{Eq:LaplaceTheorem}
	\Delta_{M} u = \Delta \hat{u}	\,.
\end{equation}
Similar to the surface gradient case, we compute a numerical approximation to $\Delta \hat{u}$ on the tangent plane. Using the same ideas as that in Section~\ref{sec:SurfaceGradient}, and specifically Eq.\,\eqref{Eq:RotationOfGradient}, we get
\begin{align}
	\Delta_{M} u &= \Delta \hat{u}\,,	\\
	&= \nabla \cdot \nabla \hat{u} \,, \\
	&= \nabla \cdot \left( R^T \nabla_{tn} \hat{u} \right)\,,\\
	&= \left( R \nabla \right) \cdot \left( R R^T \nabla_{tn} \hat{u} \right)\,, \label{Eq:DivRotInv1}\\
	&= \nabla_{tn} \cdot \left( R R^T \nabla_{tn} \hat{u} \right)\,,\\	
	&= \Delta_{tn} \hat{u}\,,\\
	&= \Delta_{T} \hat{u}\,,
\end{align}
where $\Delta_{T}$ is the $2$ dimensional Laplacian on the tangent plane, and $RR^T = R^TR = I$. Eq.\,\eqref{Eq:DivRotInv1} uses the fact that the divergence is invariant under rotations. The last equation arises due to the $\vec{n}$ direction derivatives being $0$. The same can alternatively be derived using the rotational invariance of the Laplacian.

Similar to before, we compute a $2$ dimensional volumetric numerical operator for the Laplacian on the tangential plane, which directly gives the discretization for the surface Laplacian. 
\begin{align}
	\sum_{j\in S_i}c_{ij_{T}}^{\Delta_{T}}m_{j_{T}} &= \Delta_{T} m (\vec{x}_i) \qquad \forall m\in\mathcal{P}_{T}\,,\label{Eq:TPL_Consistency}\\
	\text{min } J_i &= \sum_{j\in S_i} \left( \frac{ c_{ij_{T}}^{\Delta_{T}} } {W_{ij_{T}} } \right)^2\,. \label{Eq:TPL_Min}
\end{align}
The numerical surface Laplacian is then given by
\begin{equation}
	\label{Eq:LapBeltNumericalDefinition}
	\widetilde{\Delta}_{M} u = \sum_{j\in S_i} c_{ij}^{\Delta_M}u_j\,,
\end{equation}
where $c_{ij}^{\Delta_M} = c_{ij}^{\Delta_T}\;$.

\subsubsection{Optimized Surface Laplacian}
~\\
As mentioned earlier, all modifications of volume based meshfree GFDMs can be directly carried over to the case of surface operators. An important issue in volumetric meshfree GFDMs is the use of optimized Laplacian stencils to improve stability for Poisson problems \cite{SeiboldThesis}. For this, instead of a direct minimization according to Eq.\,\eqref{Eq:TPL_Min}, it is desired that the central stencil value $| c_{ii}^{\Delta} |$ is made ``as large as possible'' in relation to the neighbouring ones $| c_{ij}^{\Delta} |, \phantom{s} j\neq i$, while maintaining the consistency conditions of Eq.\,\eqref{Eq:TPL_Consistency}. Several procedures for the same have been done in the past (for example, \cite{Froese2018,Seibold2008}). Here, we follow the method done in our earlier work \cite[Section~2.5.5]{Suchde2018_Thesis}. A short explanation of the same is given in Appendix~\ref{App:Laplace}.

\subsection{Anisotropic Surface Laplacian}

We generalize the previous section to an anisotropic surface Laplacian~(surface diffusion) operator $\nabla_{M} \cdot \kappa \nabla_{M}\;$, for diffusion coefficeint $\kappa$. Proceeding in the same manner as earlier, the surface diffusion of a scalar valued function can be evaluated as follows
\begin{align}
	\nabla_{M} \cdot \left( \kappa \nabla_{M} u \right) &= \nabla\cdot\left( \hat{\kappa} \nabla \hat{u} \right)\,,\\
    & = \nabla \cdot \left(\hat{\kappa} R^T \nabla_{tn} \hat{u} \right)	\,,\\
    & = \left( R\nabla \right) \cdot \left( R\hat{\kappa} R^T \nabla_{tn} \hat{u} \right)	\,,\label{Eq:DivRotInv2}\\
    &= \nabla_{tn} \cdot \left( \underbrace{R\hat{\kappa} R^T}_{\hat{\kappa}_R  } \nabla_{tn} \hat{u} \right)\,,\\	
	&= \nabla_{T} \cdot \left( \hat{\kappa}_{RT}  \nabla_{T} \hat{u} \right) \,,	
\end{align}
where $\nabla_{T} \cdot \hat{\kappa}_{RT}  \nabla_{T}$ is the $2$ dimensional diffusion operator on the tangential plane, $\hat{\kappa}$ is the (possibly component-wise) normal extension of $\kappa$, and $\hat{\kappa}_{RT}$ is the appropriate submatrix of $\hat{\kappa}_{R}$. For scalar valued $\kappa$, $\hat{\kappa}_R = \hat{\kappa} = \hat{\kappa}_{RT}$. Eq.\,\eqref{Eq:DivRotInv2} uses the fact that the divergence is invariant under rotations. The last equation arises due to the $\vec{n}$ direction derivatives being $0$.

Thus, similar to the case in the earlier sections, the surface diffusion operator at a point $i$ has been reduced to the $2$ dimensional volumetric diffusion operator on the tangent plane of that point, which is approximated numerically. 
\begin{align}
	\sum_{j\in S_i}c_{ij_{T}}^{D_{T}}m_{j_{T}} &= D_{T} m (\vec{x}_i) \qquad \forall m\in\mathcal{P}_{T}\,,\label{Eq:TPD_Consistency}\\
	\text{min } J_i &= \sum_{j\in S_i} \left( \frac{ c_{ij_{T}}^{D_{T}} } {W_{ij_{T}} } \right)^2\,, \label{Eq:TPD_Min}
\end{align}
where $D_{T}$ is used as a shorthand for $\nabla_{T} \cdot  \hat{\kappa}_{RT}  \nabla_{T}\;$ . The numerical surface diffusion operator is then given by
\begin{equation}
	\label{Eq:DiffusionNumericalDefinition}
	\nabla_{M} \cdot \left( \kappa \nabla_{M} u \right) \approx 
	\widetilde{D}_{M} u = \sum_{j\in S_i} c_{ij}^{D_M}u_j\,,
\end{equation}
where $c_{ij}^{D_M} = c_{ij}^{D_T} \;$.

\subsubsection{Anisotropic Surface Laplacian with Large Jumps in Diffusion Coefficient}
\label{sec:DiffusionLargeJumps}
~\\
It is often required to model diffusion with not just a discontinuous diffusion coefficient, but one with large jumps~(with several orders of magnitude). To do the same, modifications needs to be made to the standard procedure of computing the numerical diffusion operator mentioned above. 

To achieve the same on a manifold, we extend our earlier volumetric work. For a scalar-valued $\kappa$, the following consistency conditions are enforced in addition to Eq.\,\eqref{Eq:TPD_Consistency}.
\begin{align}
	\sum_{j\in S_i}c_{ij_{T}}^{D_{T}} \frac{1 }{\kappa_{j_{T} }} &= -\Delta_{M} \left( \log \kappa \right) \,,\label{Eq:DiffExtra1}\\
	\sum_{j\in S_i}c_{ij_{T}}^{D_{T}} \frac{\delta s_{ij_{T}} }{\kappa_{j_{T} }} &= - \frac{\partial}{\partial s_{M} } \left( \log\kappa \right) \,,\label{Eq:DiffExtra2}\\
		\sum_{j\in S_i}c_{ij_{T}}^{D_{T}} \frac{\left( \delta  s_{ij_{T}} \right)^2}{\kappa_{j_{T} }} &= 2 \,,\label{Eq:DiffExtra3}	
\end{align}
where $\delta s_{ij_{T}}$ are distances along the direction of greatest change in $\kappa$, $\vec{s} = \frac{ \nabla_{M} \kappa}{ \| \nabla_{M} \kappa \|}$; and $\frac{\partial}{\partial s_{M} }$ denotes the directional surface derivative along $\vec{s}$, which is obtained numerically by rotating $\nabla_{M}\;$. The proof of the validity of the these extra conditions is given in Appendix~\ref{App:Diffusion}.

We note that this is similar in concept to the work of Yoon and Song \cite{Yoon2014} who add step functions, wedge functions and scissor functions to the polynomial test functions for GFDMs.

\section{Boundary Conditions}
\label{sec:BoundaryConditions}

One of the biggest advantages of meshfree GFDMs over particle-based meshfree methods such as SPH is the ease of handling a vast variety of boundary conditions \cite{Suchde2018_Thesis}. This is carried over to meshfree GFDMs on manifolds as well. To illustrate the same, we consider a surface Poisson equation
\begin{align}
	\Delta_M u &= f \qquad \text{in } M\,,\\
	u & = g \qquad \text{on } \partial M_1\,,\\
	\vec{\nu}\cdot\nabla_M u &= l \qquad \text{on } \partial M_2\,,
\end{align}
where $\partial M_1$ and $\partial M_2$ are parts of the manifold boundary. The discretized linear system would simply be given by
\begin{align}
	\sum_{j \in S_i} c_{ij}^{\Delta_M} u_j &= f_i \qquad \text{for } i \in M \setminus \partial M  \,,\\
	u_i &= g_i  \qquad \text{for } i \in \partial M_1 \,,\\
	\sum_{j \in S_i} c_{ij}^{M,\nu} u_j &= l_i \qquad \text{for } i \in  \partial M_2\,,
\end{align}
where $c_{ij}^{M, \nu} = \vec{\nu} \cdot \left( c_{ij}^{M,x}, c_{ij}^{M,y}, c_{ij}^{M,z}     \right)^T$. 

We note that the neighbourhoods $S_i$ for boundary points $i$ will always be ``one-sided". Similar to the case of volumetric meshfree GFDMs, the computation of differential operators on boundary points follows the exact same procedure as that for interior points. The use of ghost or virtual nodes outside the domain is \textit{not} done. More details about implementing different kinds of boundary condition using meshfree GFDMs can be found in our earlier work \cite{Suchde2018_Thesis}.

\section{Higher Dimensions and Co-Dimensions}
\label{sec:HighD_coD}

The ideas presented in this paper can easily be extended to manifolds in higher dimensions, or higher co-dimensions. The main difference is the change of the projection operator. For a $k$ dimensional manifold in $\mathbb{R}^n$, we have

\begin{equation}
 	\mathbf{P} = \mathbf{I} - \sum_{r=1}^{n-k}\vec{n}_r\vec{n}_r^T\,,
\end{equation}
where the normal space is spanned by the unit normals $\vec{n}_r$. At the discrete level, once again, the only change is the way of projecting neighbouring points to the tangent space. Distances in the tangential space can still be computed by rotating actual distances in the embedding space. The only difference would be in the rotation matrix, which would computed as
\begin{equation}
	R^T = \left(\begin{array}{cccccc} 
	\vec{t}_1 & \cdots & \vec{t}_k & \vec{n}_1 & \cdots & \vec{n}_{n-k} 
	\end{array}\right)\,,
\end{equation}
where $\vec{t}_1, \dots ,\vec{t}_k$ are orthogonal to each other and each $\vec{n}_r$. Thus, $R^T$ has mutually orthogonal columns. Once the distances in the tangential space are known~(given by the appropriate submatrix of $R\, \delta\vec{x}_{ij} $), the differential operators can be computed in the same manner as volumetric $k$-dimensional operators.

\section{Numerical Results and Validation}
\label{sec:NumResults}

We emphasize that we do not use uniformly spaced point clouds, with the exception of the first example. All irregularly spaced point clouds are setup in a manner similar to that done in several volumetric meshfree GFDMs. Starting from a CAD file for the geometry, points are placed using an advancing front technique for point clouds, like that done by Drumm \textit{et~al.} \cite{Drumm2008}. Using the distance conventions used in \cite{Drumm2008, Suchde2017_CCC} results in about $15-20$ points in each neighbourhood. Normal computation follows the procedure used for free surfaces of volumetric flow by meshfree GFDM, as done by \cite{Edgar2017}.

Unless specified otherwise, in all numerical examples monomials up to the second order are considered in the computation of all numerical differential operators. Further, the sparse linear systems arising in each example are solved with a BiCGSTAB iterative solver \cite{BiCGSTAB} without the use of any pre-conditioner. If available, the solution of the previous time level is used as an initial guess.

\subsection{Heat Equation on a Sphere}

As a validation case, we consider the surface diffusion equation 
\begin{equation}
	\label{Eq:HeatEquation}
	u_t = \Delta_{M} u \,,
\end{equation}
on a unit sphere. For initial conditions $u_0 = xy$, the analytical solution to Eq.\,\eqref{Eq:HeatEquation} is given by
\begin{equation}
	\label{Eq:AnalSurfHeat}
	u_{exact}(\vec{x},t) = \exp(-6t)xy \,.
\end{equation}
We note that the same example has also been considered in Chen \textit{et~al.} \cite{Chen2014}. The sphere is discretized with a quasi-uniformly distributed point cloud generated by DistMesh \cite{Persson2004}. Irregularly spaced point clouds are considered in the later sections. For consistency with the case of irregularly spaced point clouds considered from the next example, here, neighbourhoods are taken to be the $15$ closest points, including the center point. Further, $h$ is taken as the maximum distance between a center point $i$, and all its neighbouring points $j\in S_i$.

The point cloud is taken to be quasi-regularly spaced in this example only to determine the numerical order of convergence. Further, a small time step of $\Delta t = 0.1 h^2$ is used for the same reason. For a numerical solution $u$, relative errors in the solution are measured at $t=0.3$ as
\begin{equation}
	\epsilon_2 = \left[ \frac{\sum_{i=1}^N \left( u_i - u_{exact}(\vec{x}_i) \right)^2}{\sum_{i=1}^N \left( u_{exact}(\vec{x}_i) \right)^2 } \right]^{\frac{1}{2}}\,.
\end{equation}

\begin{figure}
  \centering
  \includegraphics[width=0.7\textwidth]{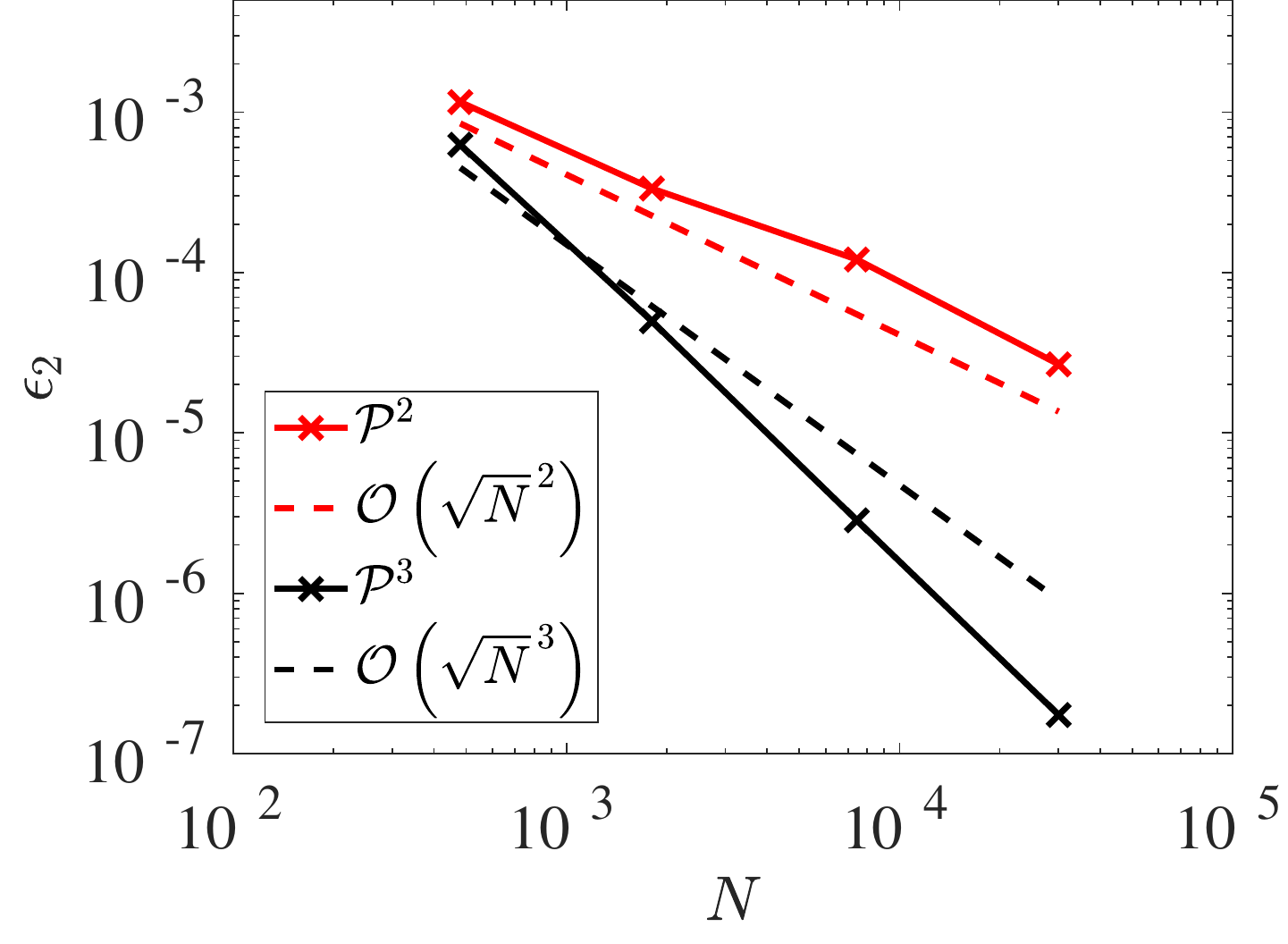}
  \caption{Surface heat equation on a sphere: Relative errors against number of points $N$ in the domain. Second and third order convergence rates are marked with dashed lines. Numerical errors are shown by solid lines. The superscript of $\mathcal{P}^p$ indicates the maximum order of monomials being used in the computation of the differential operators.}%
  % Figure created in SurfaceDiffusionEquationWrapper.m
  \label{Fig:SurfaceHeatEquation}
\end{figure}

A Crank--Nicolson time integration method is used. We consider two cases here: one with monomials up to the second order, and one with monomials up to the third order. The plots of the relative error against the number of points $N$ are shown in Figure~\ref{Fig:SurfaceHeatEquation}. In each of the cases, the experimental order of convergence matches the theoretical expectation. Due to the use of small time steps, the BiCGSTAB solver converged within very few iterations. Up to a tolerance of $10^{-10}$, the maximum number of iterations required for convergence was $5$, $5$, $4$ and $3$ for the second order case; and $3$, $3$, $3$ and $2$ for the third order case. In both cases, the number of iterations are reported in the order of increasing $N$. % Averages: 2.67, 1.66, 1.05, 0.54; and 1.714, 0.67, 0.82, 0.53  (Note: MATLAB count also gives half iterations.)

To obtain an order of accuracy higher than third order, the size of the neighbourhoods would need to be increased. To use monomials up to order $p$, for a $k$ dimensional manifold, the number of monomial functions needed is given by the binomial coefficient ${p+k}\choose{k}$. Thus, to ensure that the local least squares systems for computing the differential operators are solvable, a necessary condition is that the number of neighbours should be more than ${p+k}\choose{k}$. The increased support sizes results in denser linear systems. This is one of the limitations of using classical GFDMs to approximate derivatives. To overcome this, modification in volumetric GFDMs to obtain higher order accuracy \cite{Milewski2012,Trask2017} would need to be extended to surfaces. 
\subsection{Diffusion on a Torus with forcing}

We consider the diffusion equation again, but with forcing
\begin{equation}
	\label{Eq:HeatEquationF}
	u_t = \Delta_{M} u + f(u, \vec{x}, t)\,.
\end{equation}
The domain is taken to be a torus given by
\begin{equation}
	\left( 1 - \sqrt{ x^2 + y^2 } \right)^2 + z^2 = \frac{1}{9}\,.
\end{equation}
A manufactured solution is considered, as done in \cite{Fuselier2013}. The exact solution is given by
\begin{equation}
\label{Eq:TorusExact}
	u_{exact}(\vec{x}, t) = \frac{1}{8}\exp{(-5t)}x\left(x^4 - 10x^2y^2 + 5y^4 \right)\left( x^2 + y^2 - 60z^2 \right)\,.
\end{equation}
The forcing function is taken such that Eq.\,\eqref{Eq:TorusExact} satisfies Eq.\,\eqref{Eq:HeatEquationF}, and is given in \cite{Fuselier2013}. Once again, a Crank--Nicolson time integration method is used with a small $\Delta t$, as done above. From this example onwards, all points clouds are taken to be irregularly spaced, and are set up in a manner similar to \cite{Drumm2008}. The plots of relative error against the number of points $N$ are shown in Figure~\ref{Fig:SurfaceHeatEquation}. The figure illustrates that projection to the tangent plane along the central normal~(as done in this paper) is more accurate in the present context than the projection along the neighbour normals~(as done in some mesh-based methods \cite{Demanet2006}). Further, the experimental order of convergence is seen to match the theoretical expectation of second order for both types of projection. Up to a tolerance of $10^{-10}$, the maximum number of iterations required for convergence was $7$, $7$, $6$ and $5$, in order of increasing $N$. We note that the larger time steps used in the coming examples require higher number of iterations for convergence of the sparse linear system solver.
% for the central normal projection case; and $7$, $6$, $4$ and $7$ for the neighbour normal projection case. In both cases, the number of iterations are reported in the order of increasing $N$. We note that the larger time steps used in the coming examples result in more number of iterations for convergence of the sparse linear system solver. % Averages: 6.28, 4.98, 3.99, 2.05; and 6.25, 4.95, 3.5, 0.67  (Note: MATLAB count also gives half iterations.)
%
\begin{figure}
  \centering
  \includegraphics[width=0.7\textwidth]{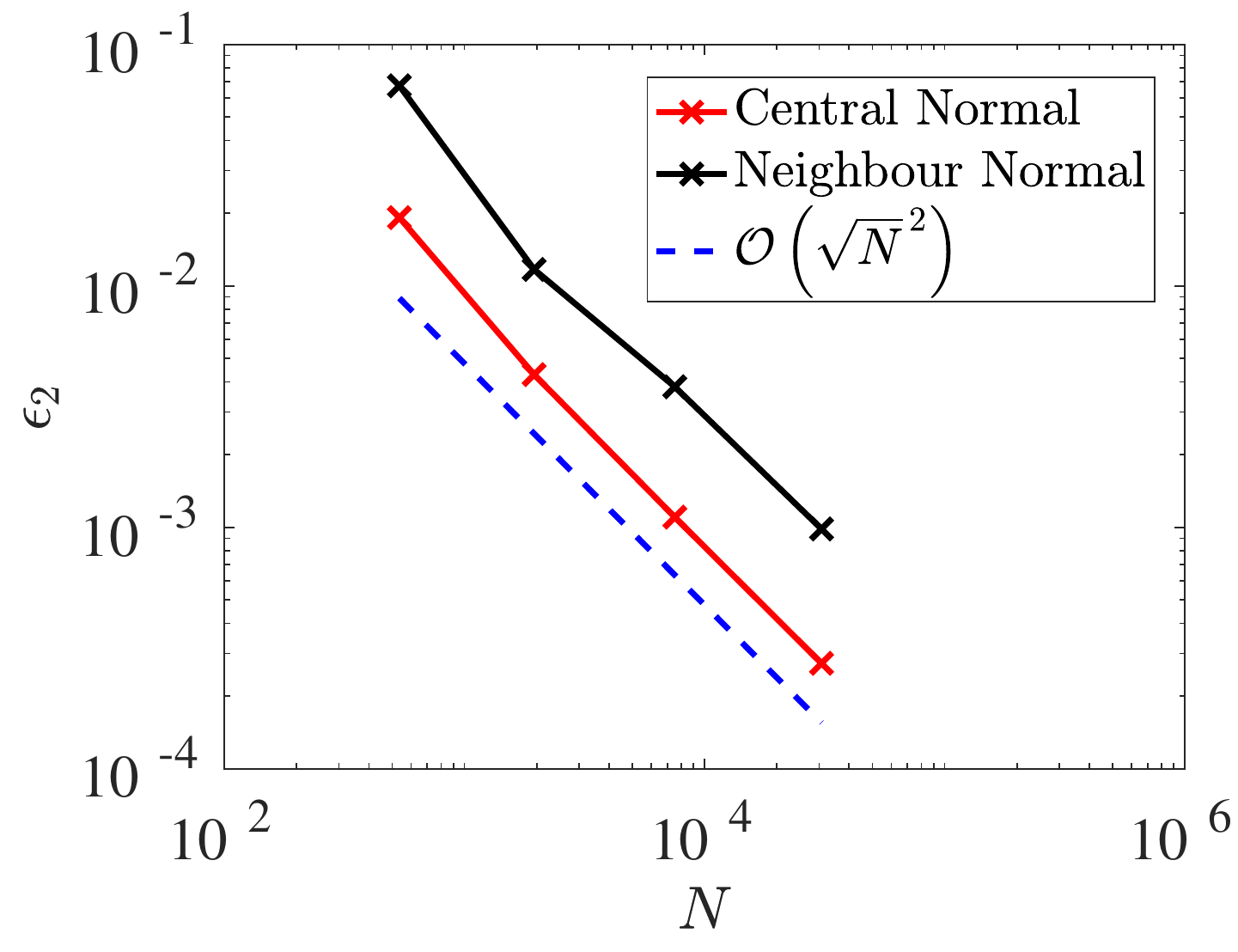}
  \caption{Diffusion on a torus with forcing: Relative errors against the number of points $N$ in the domain. Second order convergence rate is marked in blue. Numerical errors are shown in red~(projection to tangential plane along the central normal, as done in this paper), and in black~(projection along neighbour normals, as done in some mesh-based work \cite{Demanet2006}).}%
  % Figure created in SurfaceDiffusionEquation2Wrapper.m
  \label{Fig:SurfaceHeatEquationTorus}
\end{figure}

\subsection{Four Strip Problem on a Surface}

We now consider an elliptic problem with rough coefficients
\begin{equation}
\label{Eq:DivEtaGradPhi}
- \nabla_{M}\cdot\eta\nabla_{M} \phi = f\,,
\end{equation}
with appropriate boundary conditions. A common case to test numerical schemes is the volumetric equivalent of Eq.\,\eqref{Eq:DivEtaGradPhi} with discontinuous $\eta$ (for example, \cite{Trask2017}). Here, we consider not just rough coefficients, but $\eta$ with large jumps, up to several orders of magnitude. 

\begin{figure}
  \centering
  \includegraphics[width=0.7\textwidth]{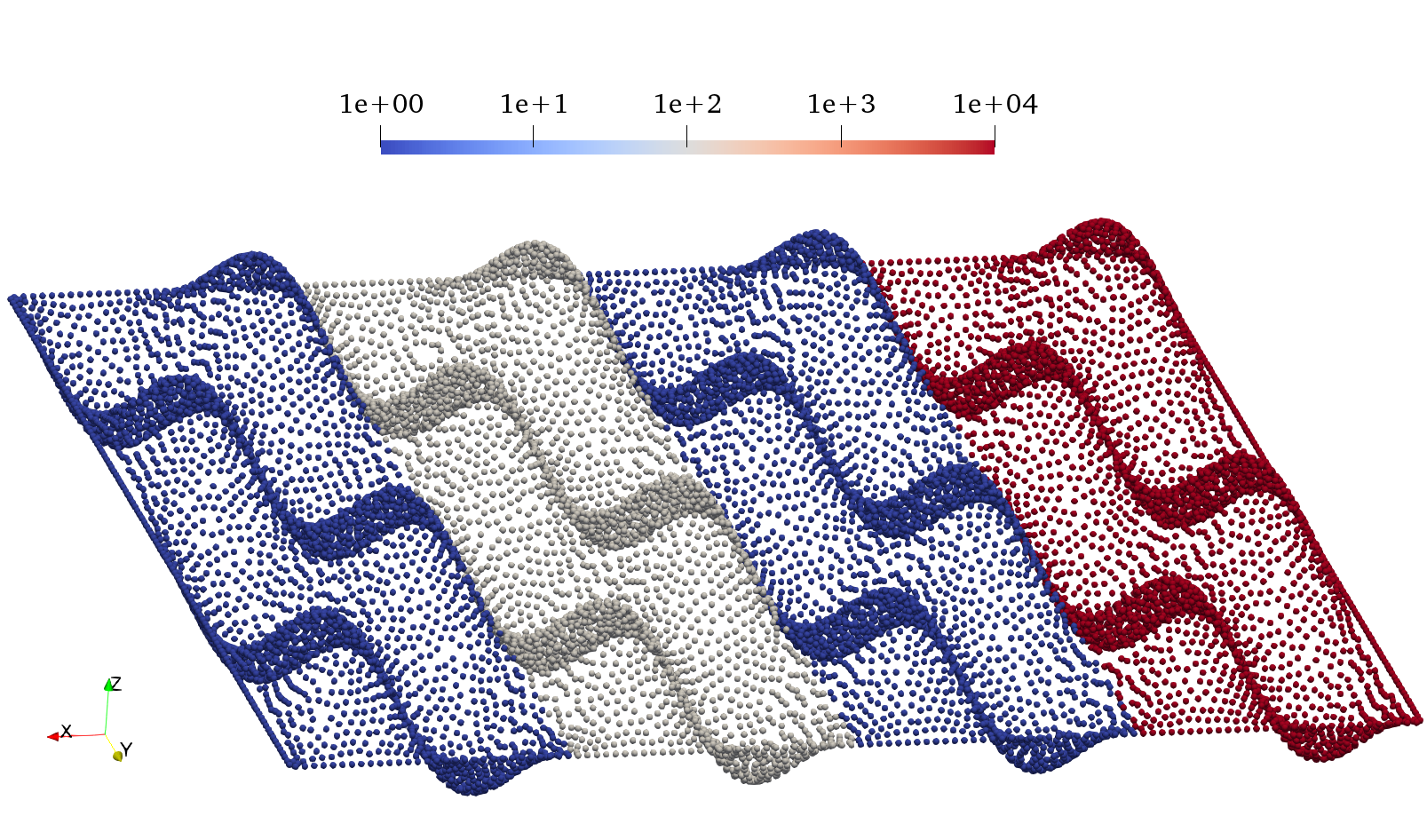}
  \caption{Discretized domain for the elliptic four-strip problem given by $z=\sin(2x)\sin(y)$. The colour indicates the different values of $\eta$.}%
  % Figure created in Paraview - see SaveState3_forDomain.pvsm in Manifold_Results/FiveStripProblem/
  \label{Fig:DivEtaGradProblem_Domain}
\end{figure}

Eq.\,\eqref{Eq:DivEtaGradPhi} with $f\equiv 0$ is solved on a wave surface defined by $z=\sin(2x)\sin(y)$, $x\in[0,4\pi]$, $y\in[0,4\pi]$, which is shown in Figure~\ref{Fig:DivEtaGradProblem_Domain}. The domain has four strips along the $x$ direction, each of the same width, with different values of the diffusion coefficient $\eta$. Starting from $x=0$, $\eta_1=10^4$, $\eta_2=10^0$, $\eta_3=10^2$, and $\eta_4 =10^0$. Figure~\ref{Fig:DivEtaGradProblem_Domain} shows the different strips for $\eta$. It also illustrates that the points are unevenly distributed, and that extra points are \textit{not} added along the interfaces where $\eta$ changes. 

\begin{figure}
  \centering
  \includegraphics[width=0.48\textwidth]{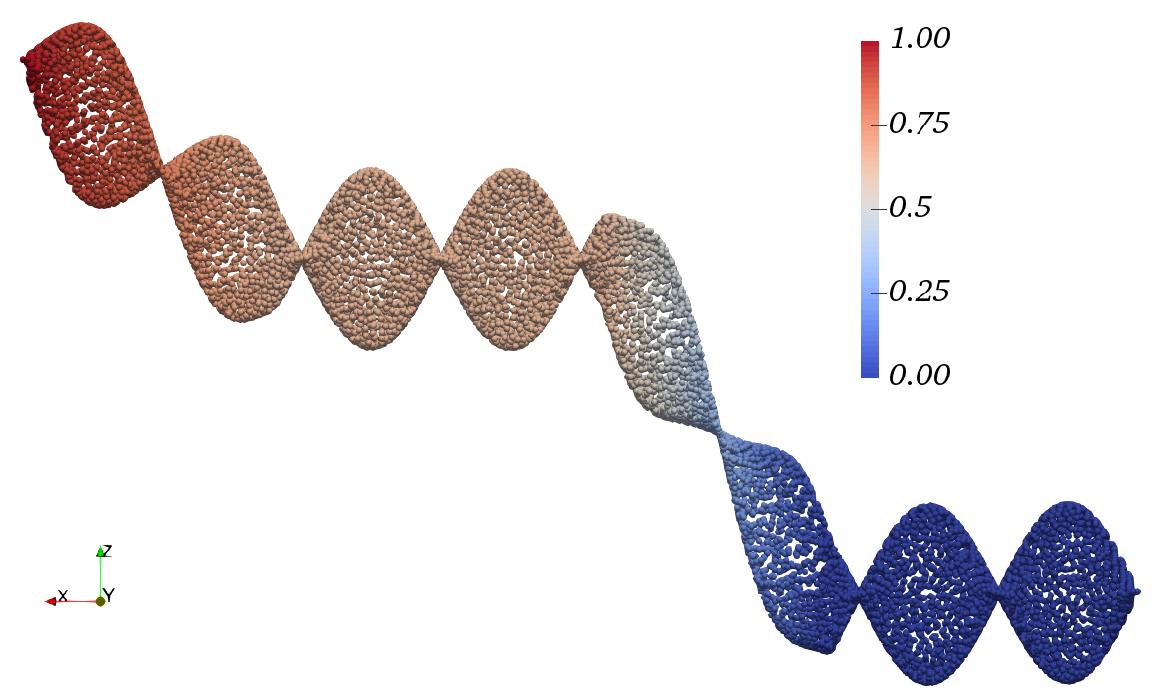}
  \includegraphics[width=0.48\textwidth]{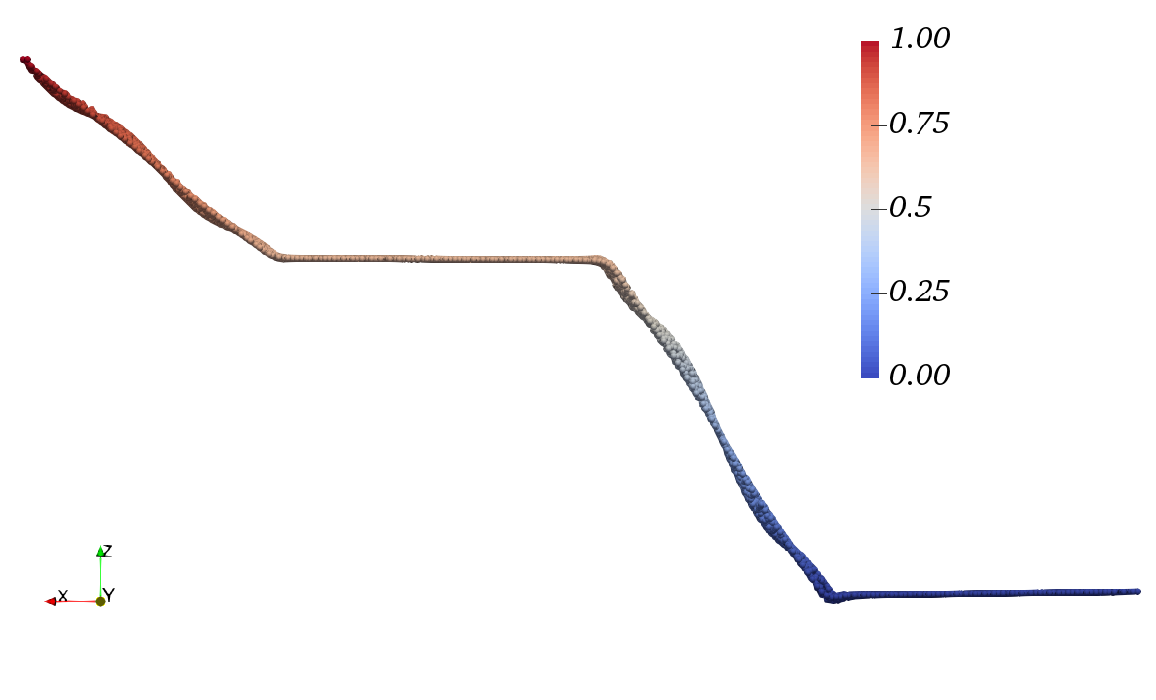}
  \caption{Results for the elliptic four-strip problem: Domain warped by the solution. $x$ vs. $y$ vs. $z+6\phi$~(left) and $x$ vs. $y$ vs. $6\phi$~(right). The colour of the points indicates the value of $\phi$.}%
  % Figures created in Paraview - see SaveState2.pvsm in Manifold_Results/FiveStripProblem/
  \label{Fig:DivEtaGradProblem_Solution}
\end{figure}

Dirichlet boundary conditions are enforced at two ends, and pseduo-Neumann boundary conditions on the other two.

\begin{align}
\phi &= 0\phantom{abc} \text{if }\, x=0\,,\\
\phi &= 1\phantom{abc} \text{if }\, x=4\pi\,,\\
\vec{\hat{y}}\cdot\nabla_{M} \phi &= 0\phantom{abc} \text{if }\, y=0\,,\\
\vec{\hat{y}}\cdot\nabla_{M} \phi &= 0\phantom{abc} \text{if }\, y=4\pi\,,
\end{align}
where $\vec{\hat{y}}=(0,1,0)$. The domain is discretized with a total of $N=13918$ points. Due to the large jumps in $\eta$, the addition of the extra test functions for the computation of the diffusion operator, as explained in Section~\ref{sec:DiffusionLargeJumps}, is essential to prevent excessive numerical oscillations. Using these, the numerical solution is plotted in Figure~\ref{Fig:DivEtaGradProblem_Solution}, where the domain is warped by a scalar multiple of the solution $\phi$. The figure illustrates that the result shows a good agreement with the expectation that the solution does not vary in the $y$ direction, and that the linear slope of the solution in the $x$ direction is inversely related to the values of $\eta$.

\subsection{Advection on a manifold}

Having considered a parabolic and an elliptic PDE, we now consider the hyperbolic problem of the transport equation
\begin{equation}
	\label{Eq:Advection}
	\frac{\partial \phi}{\partial t} + \vec{v}\cdot \nabla_{M}\phi = 0\,,\
\end{equation}
on the surface of a cone given by
\begin{equation}
 x^2 + y^2 = \frac{4}{9}z^2\,,
\end{equation}
with $-6 \leq z \leq 0$. To discretize Eq.\,\eqref{Eq:Advection}, we once again take advantage of the fact that developments in volumetric GFDMs can be directly carried over here. Numerical methods to discretize the advection term for the volumetric case have been widely studied, including using meshfree GFDMs \cite{Praveen2007, SeifarthThesis}. Here, we follow the work of Seifarth \cite{SeifarthThesis} for volumetric meshfree GFDMs. A upwind discretization of Eq.\,\eqref{Eq:Advection} leads to the following semi-discrete form
\begin{equation}
	\frac{d\phi_i}{dt} = -2 \sum_{\substack{j\in S_i \\ j \neq i}} c_{ij}^{\vec{v}\cdot\nabla_M} (\phi_{ij} - \phi_i) \,,
\end{equation}
where
\begin{align}
	c_{ij}^{\vec{v}\cdot\nabla_M}  &= c_{ij}^{M,x}v_i^1 + c_{ij}^{M,y}v_i^2 + c_{ij}^{M,z}v_i^3\,,\\
	\phi_{ij} &= \frac{1}{2} \left[ \left( 1 + \text{sign}(\delta\vec{x}_{ij} \cdot \vec{v}_i ) \right)\phi_{ij}^{+}  +  \left( 1 - \text{sign}(\delta\vec{x}_{ij} \cdot \vec{v}_i ) \right)\phi_{ij}^{-} \right]\,,
\end{align}
where $\delta\vec{x}_{ij} = \vec{x}_j - \vec{x}_i$ are the distances computed in $\mathbb{{R}}^3$; and $\phi_{ij}^{+}$ and $\phi_{ij}^{-}$ are the reconstructed values. Here, we use a MUSCL reconstruction  with a Superbee limiter. A pure upwind scheme is also shown for comparison. For both cases, a SDIRK2 \cite{Alexander1977, Ellsiepen1999} implicit second order method of time integration is used, which enables the use of large time steps. A linearization of the arising sparse implicit system is done before using a BiCGSTAB solver. More details of the scheme, including the time integration, can be found in \cite[Section 4.3]{SeifarthThesis}. % See page 109 of Seifarth thesis.

%The cone domain is of height $6$ and base radius $4$. The vertex of the cone is at the origin. 
A Gaussian bell is transported on the cone. The initial condition is given by
\begin{equation}
	\label{Eq:Advection_IC}
	\phi(\vec{x},0) = \left\{ \begin{array}{ll}
         \frac{ \exp\left( - \| \vec{x} - \vec{x}_0\|^2\right) - \exp\left(-25\right) }{1 - \exp\left( -25 \right) } &\; \text{if } \| \vec{x} - \vec{x}_0\|^2 < 5^2\,, \\
          0,      &\; \text{elsewhere }\,,
        \end{array} \right.
\end{equation}
for $ \vec{x}_0 = (2,0,-3)$. The velocity field is taken to be $\vec{v} = (-y,x,0)$. Note that this is not divergence-free on the manifold, and that it lies on the tangent bundle of the manifold.

\begin{figure}
  \centering
  \includegraphics[width=0.3\textwidth]{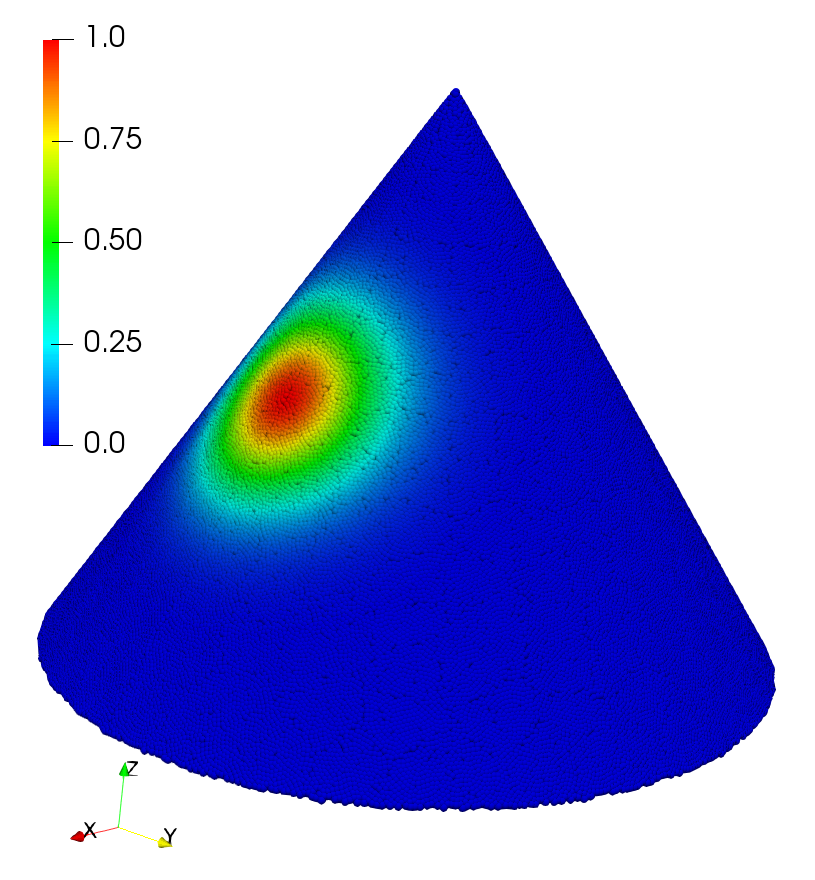}\phantom{abcdefghij}
  \includegraphics[width=0.3\textwidth]{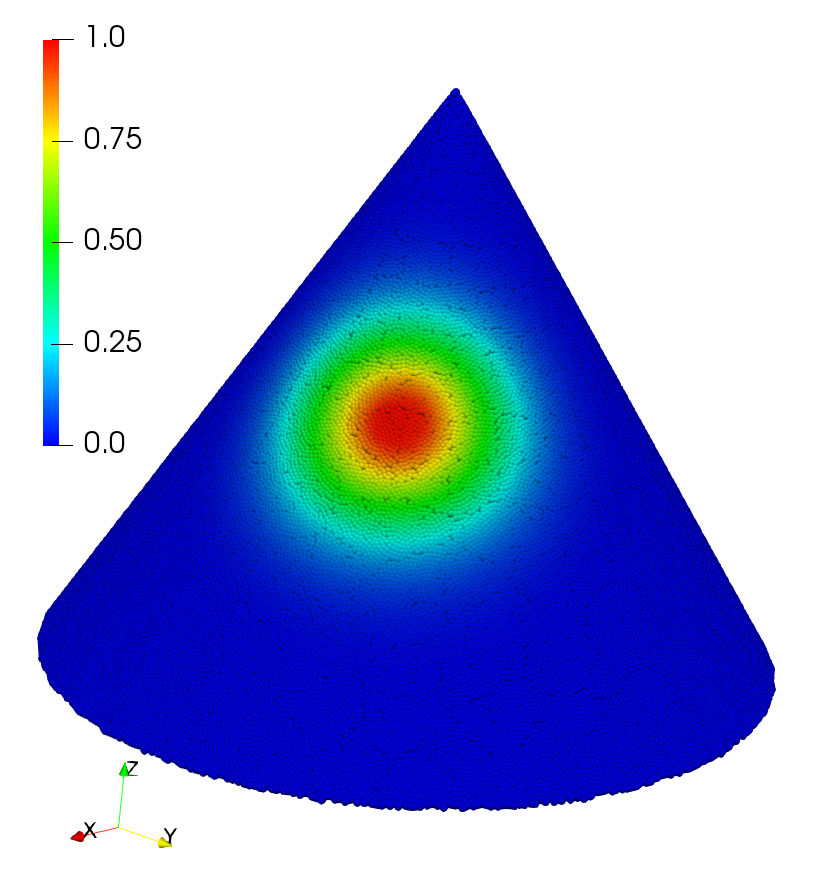}\\
  \includegraphics[width=0.3\textwidth]{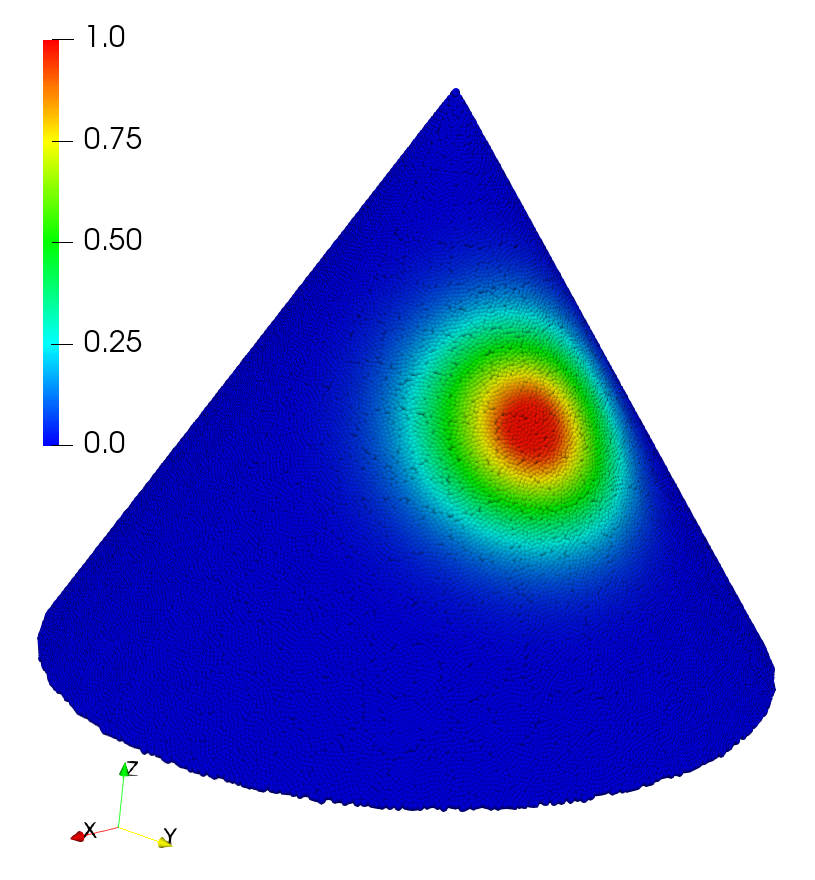}\phantom{abcdefghij}
  \includegraphics[width=0.3\textwidth]{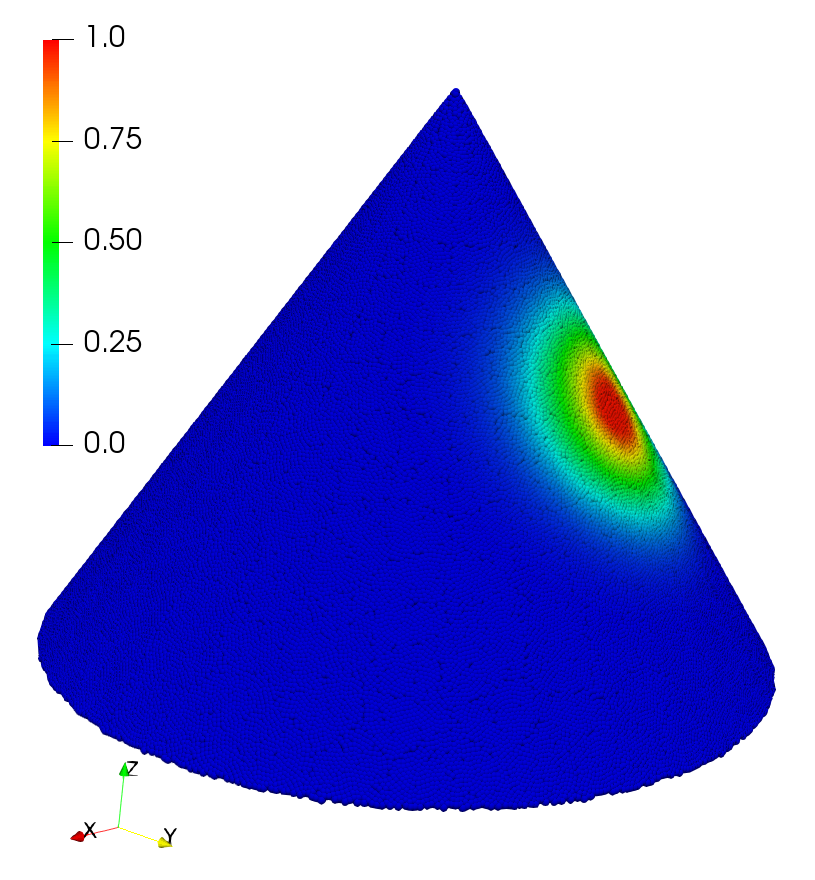}\\  
  \caption{Advection on the surface of a cone for $N=52\,193$: evolution of the solution $\phi$ at different times using a MUSCL reconstruction with a Superbee limiter. $t=0 s$~(top left), $t=0.66 s$~(top right), $t=1.32 s$~(bottom left), and $t=1.98 s$~(bottom right).}%
  % Figures created in Paraview, from DISS_Suchde/Manifold_Results/ManifoldTransport/NewResults_ConeHeight6/Ht6_h0.115_Limiter3
  \label{Fig:Advection_Results}
\end{figure}

\begin{figure}
  \centering
  \includegraphics[width=0.3\textwidth]{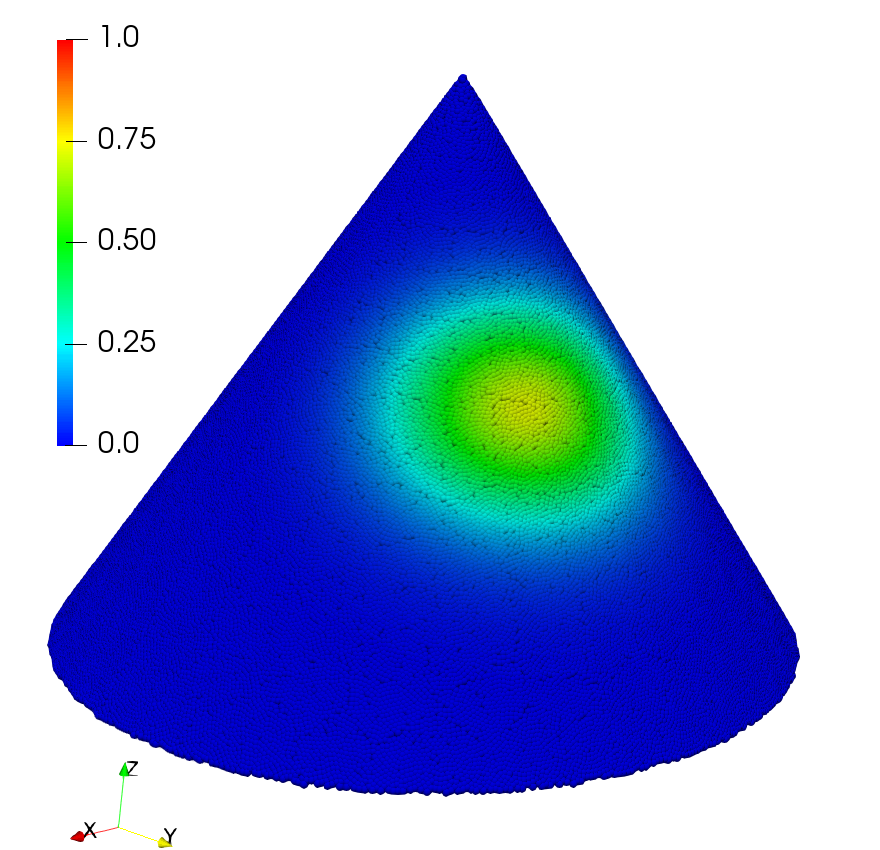}\phantom{abcdefghij}
  \includegraphics[width=0.3\textwidth]{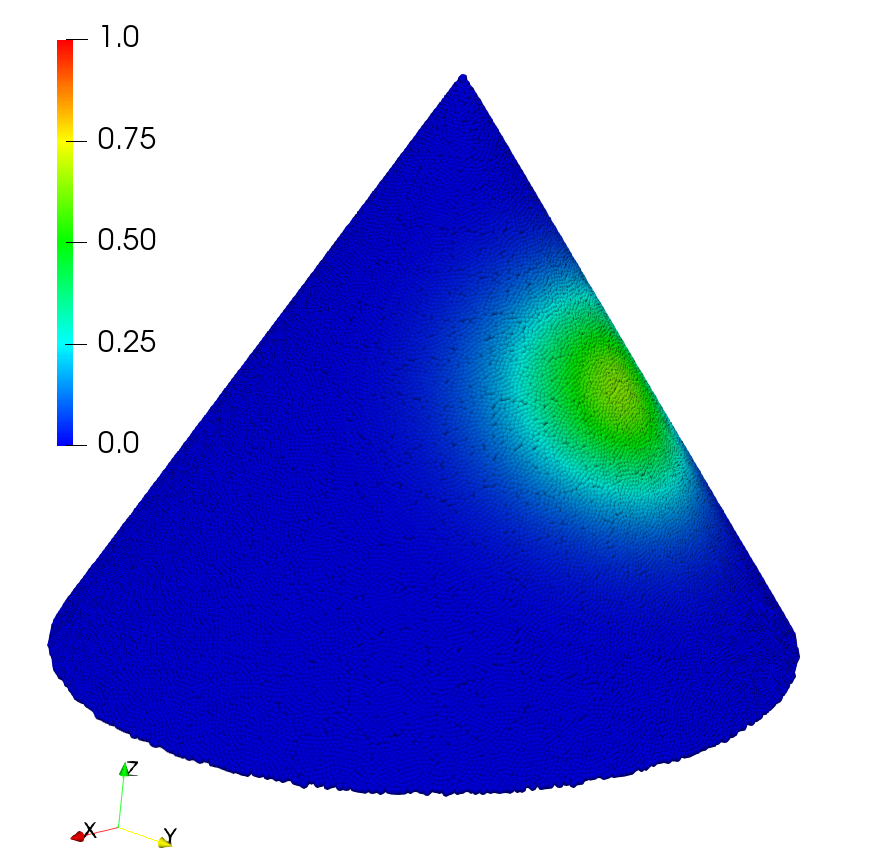} 
  \caption{Advection on the surface of a cone for $N=52\,193$: evolution of the solution $\phi$ at different times for a pure upwind scheme. $t=1.32 s$~(left), and $t=1.98 s$~(right).}%
  % Figures created in Paraview, from DISS_Suchde/Manifold_Results/ManifoldTransport/NewResults_ConeHeight6/Ht6_h0.115_Limiter0
  \label{Fig:Advection_PureUpwind}
\end{figure}

Time integration is performed with a time step of $\Delta t = 0.02$ until $t_{end} = 2\pi$, which corresponds to one full rotation. The domain is discretized with $N=52\,193$ points. The results using the scheme mentioned above are shown in Figure~\ref{Fig:Advection_Results}. Results for a pure upwind scheme~(without any reconstruction) with the same time integration method and same initial conditions are shown in Figure~\ref{Fig:Advection_PureUpwind}. These figures illustrate that the MUSCL reconstruction with Superbee limiter results in minimal numerical diffusion, while the pure upwind scheme causes excessive numerical diffusion. 
Further, a direct central difference approximation of the $\vec{v}\cdot\nabla$ operator leads to unstable simulations. Each of these observations agrees with the expectation for volumetric advection. The errors can be quantified as follows. For the numerical domain with $N=52\,193$, as used in Figures~\ref{Fig:Advection_Results} and \ref{Fig:Advection_PureUpwind}, at $t=0$, $\text{max}(\phi) = 0.99$ for both simulations (the maximum value of $1$ of Eq.\,\eqref{Eq:Advection_IC} is not attained as no point is present exactly at $\vec{x}_0$). At $t=2\pi$, $\text{max}(\phi) = 0.36$  for the pure upwind case, and $\text{max}(\phi) = 0.97$ for the MUSCL and Superbee case. %An analytical solution to the above problem can be seen to be
%
%
%\begin{equation}
%	\label{Eq:Advection_Solution}
%	\phi(\vec{x},t) = \left\{ \begin{array}{ll}
%         \frac{ \exp\left( - \| \vec{x} - \vec{x}_t\|^2\right) - \exp\left(-25\right) }{1 - \exp\left( -25 \right) } &\; \text{if } \| \vec{x} - \vec{x}_t\|^2 < 5^2\,, \\
%          0,      &\; \text{elsewhere }\,,
%        \end{array} \right.
%\end{equation}
%
%with $\vec{x}_t=(2\cos(t), 2\sin(t), -3)$. The evolution of the relative errors with the analytical solution are plotted in Figure~\ref{Fig:AdvectionErrors} for both the pure upwind and MUSCL Superbee case. A sharp increase in error is observed in the first few time steps in both cases due to numerical smoothing around the base of the Gaussian bell $\|\vec{x}-\vec{x}_t \| = 5$.
%
%

\begin{figure}
  \centering
  \includegraphics[width=0.47\textwidth]{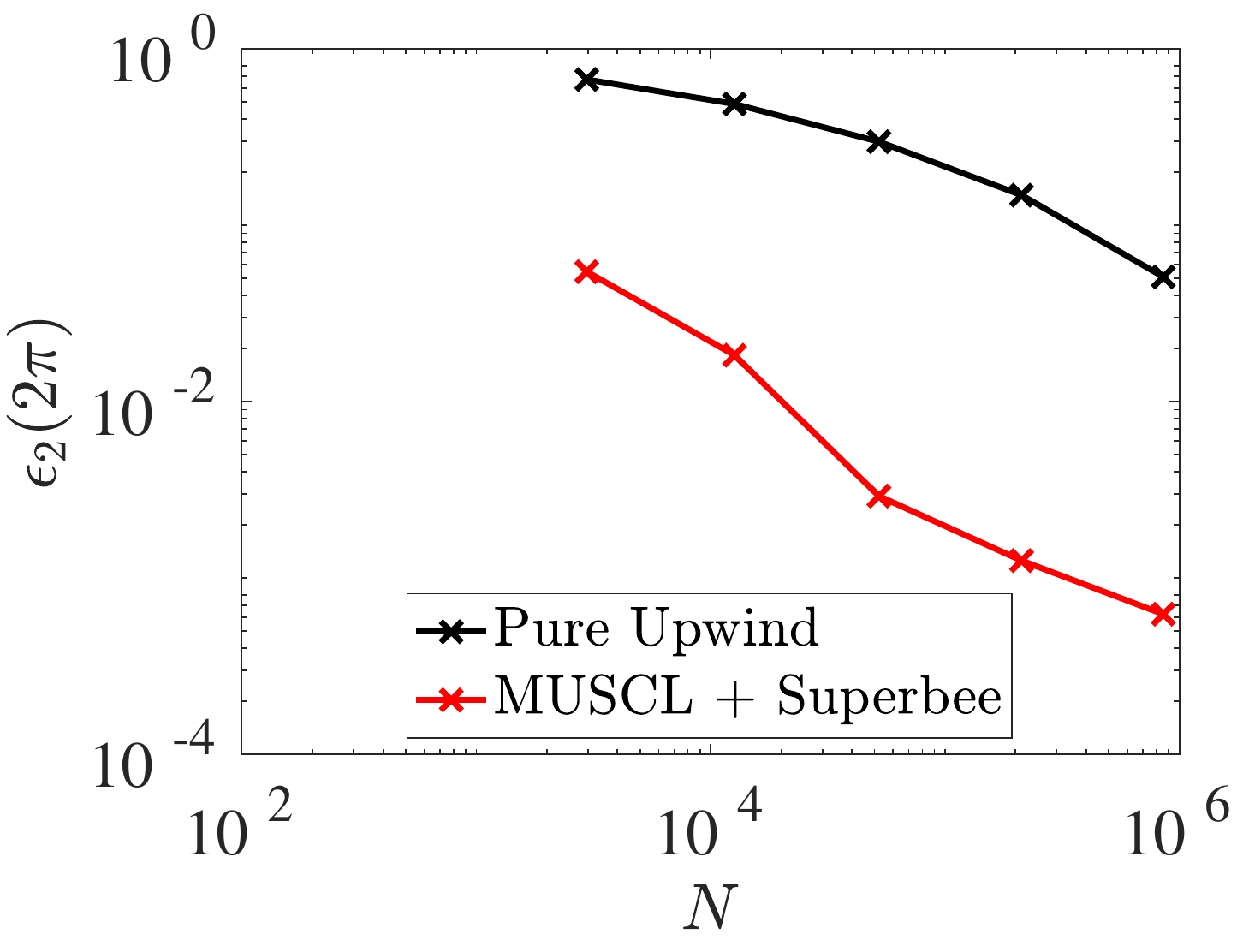}
  \includegraphics[width=0.47\textwidth]{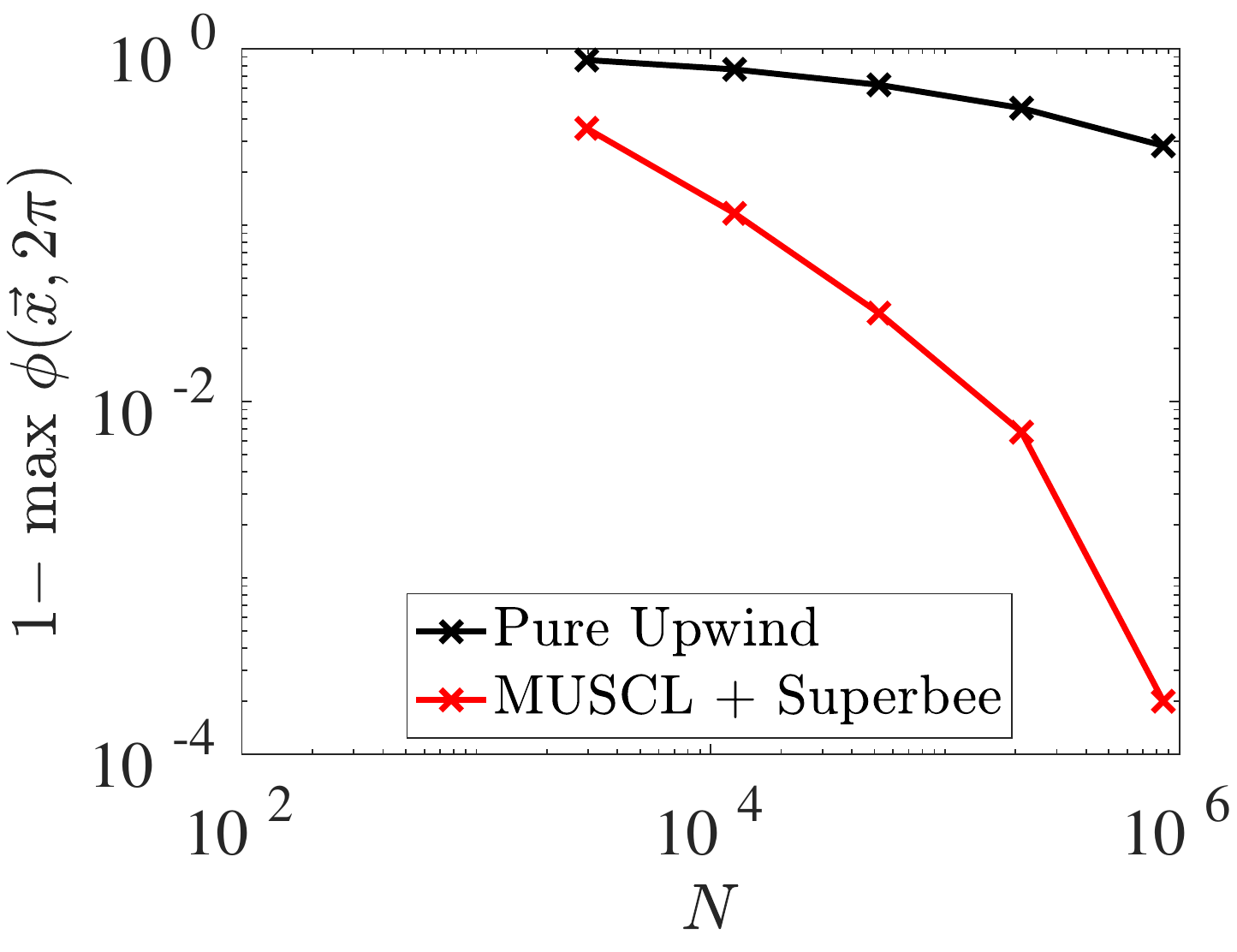}
  \caption{Advection on the surface of a cone: Relative errors~(left) and the error in the peak of the solutions~(right) at $t=2\pi$ for varying $N$. For the pure upwind case~(black), and the MUSCL + Superbee case~(red).}%
  % Figure created in SurfaceDiffusionEquation2Wrapper.m
  \label{Fig:AdvectionErrors}
\end{figure}

After, one full rotation, at $t=2\pi$, the analytical solution matches the initial condition. We measure relative errors for this with varying $N$. The plots of the same are shown in Figure~\ref{Fig:AdvectionErrors}. A smaller time step is used in these simulations $\Delta t = 0.1h$. The smoothing lengths considered in the simulations for Figure~\ref{Fig:AdvectionErrors} start with $h=0.6$ and are consecutively halved till $h=0.075$, which results in the number of points varying from $N=2960$ to $N=853837$. The value of the peak of the Gaussian is also shown in the same figure to illustrate the numerical diffusion of the peak of the solution.

This example, once again, illustrates the ease of taking over volume based GFDM ideas to the surface PDE discretization setting explained in this paper. In fact, for the simulations in this section, not only was the work of volumetric methods carried over, the code of Seifarth \cite{SeifarthThesis} was also used directly. 

\subsection{Cahn--Hilliard Equation}

We now consider the PDE system of the Cahn-Hilliard equation \cite{Cahn1958} on a closed manifold (with no boundaries). It describes the process of phase separation \cite{Kim2016}. Numerical solutions to the equation on manifolds has had a lot of interest \cite{Chen2014, Du2011}. Here, we consider the Cahn-Hilliard equation as a pair of coupled second order PDEs, as done in \cite{Gera2017}.
\begin{align}
	\frac{\partial f}{\partial t} &= \frac{1}{Pe} \Delta_{M} \mu \label{Eq:CahnHilliard1}\,,\\
	\mu &= \frac{dg}{df} - Cn^2 \Delta_{M}f\label{Eq:CahnHilliard2}\,,
\end{align}
where $f$ is an order parameter on the manifold, $\mu$ is the chemical potential, $Pe$ is the surface Peclet number, and $Cn$ is the Cahn number. As commonly done, the function $g$ is taken according to a double well potential $g(f) = \frac{f^4}{4} - \frac{f^2}{2}$. 

\begin{figure}
  \centering
  \includegraphics[width=0.5\textwidth]{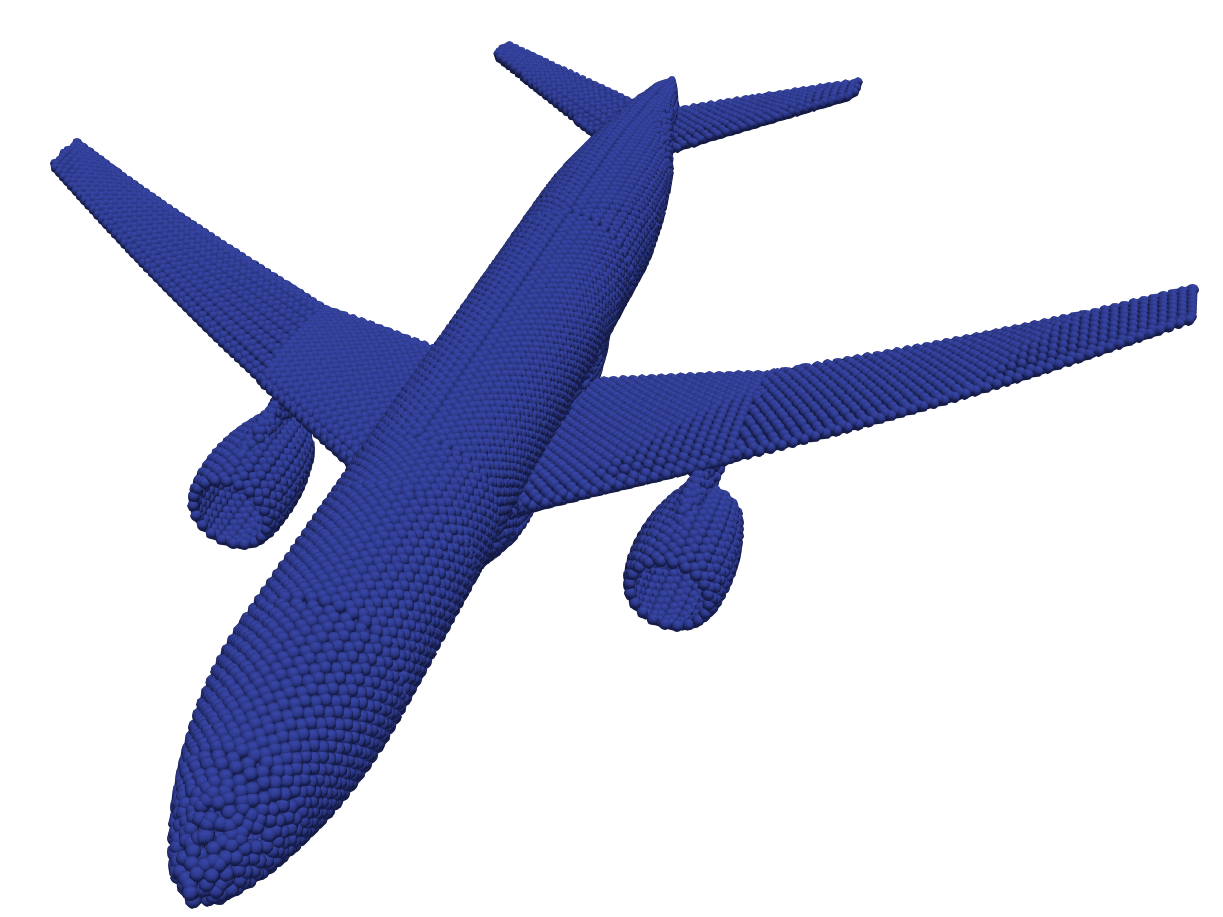}
  \caption{Discretized domain for the Cahn-Hilliard equation.}%
  % Figure created in Paraview - see SaveState3_forDomain.pvsm in Manifold_Results/FiveStripProblem/
  \label{Fig:Airplane_Domain}
\end{figure}

Eqs.~\eqref{Eq:CahnHilliard1} and \eqref{Eq:CahnHilliard2} are solved on a surface of an airplane, as shown in Figure~\ref{Fig:Airplane_Domain}. The airplane CAD files are based on that by NASA CRM \cite{Vassberg2008}. Note that each part of the airplane, including the wings, have non-zero thickness, and there is a layer of points on each side of the geometry. Thus, the discretized airplane considered here has no boundaries, and forms a closed manifold. 

\begin{figure}
  \centering
  \includegraphics[width=0.48\textwidth]{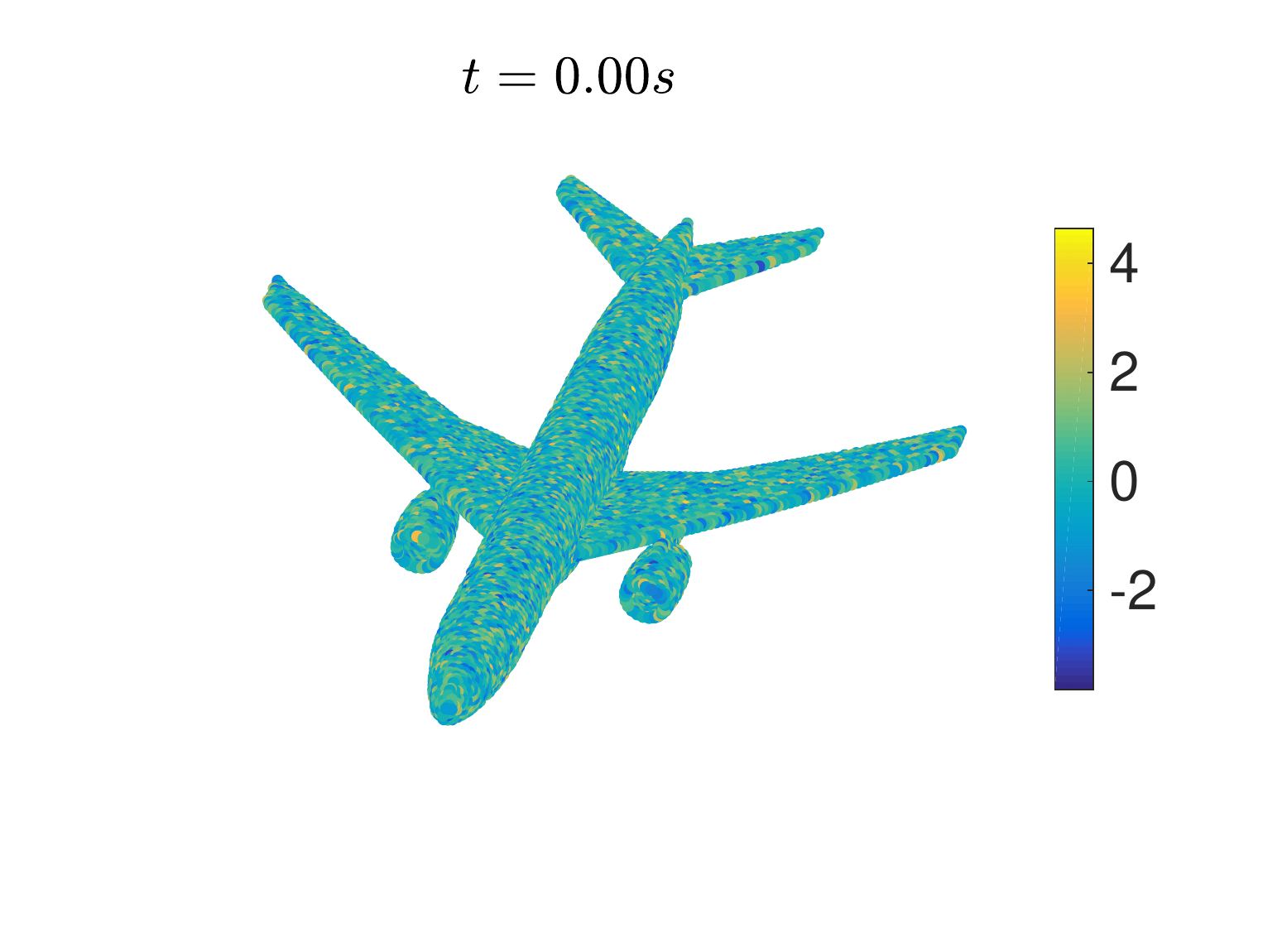}
  \includegraphics[width=0.48\textwidth]{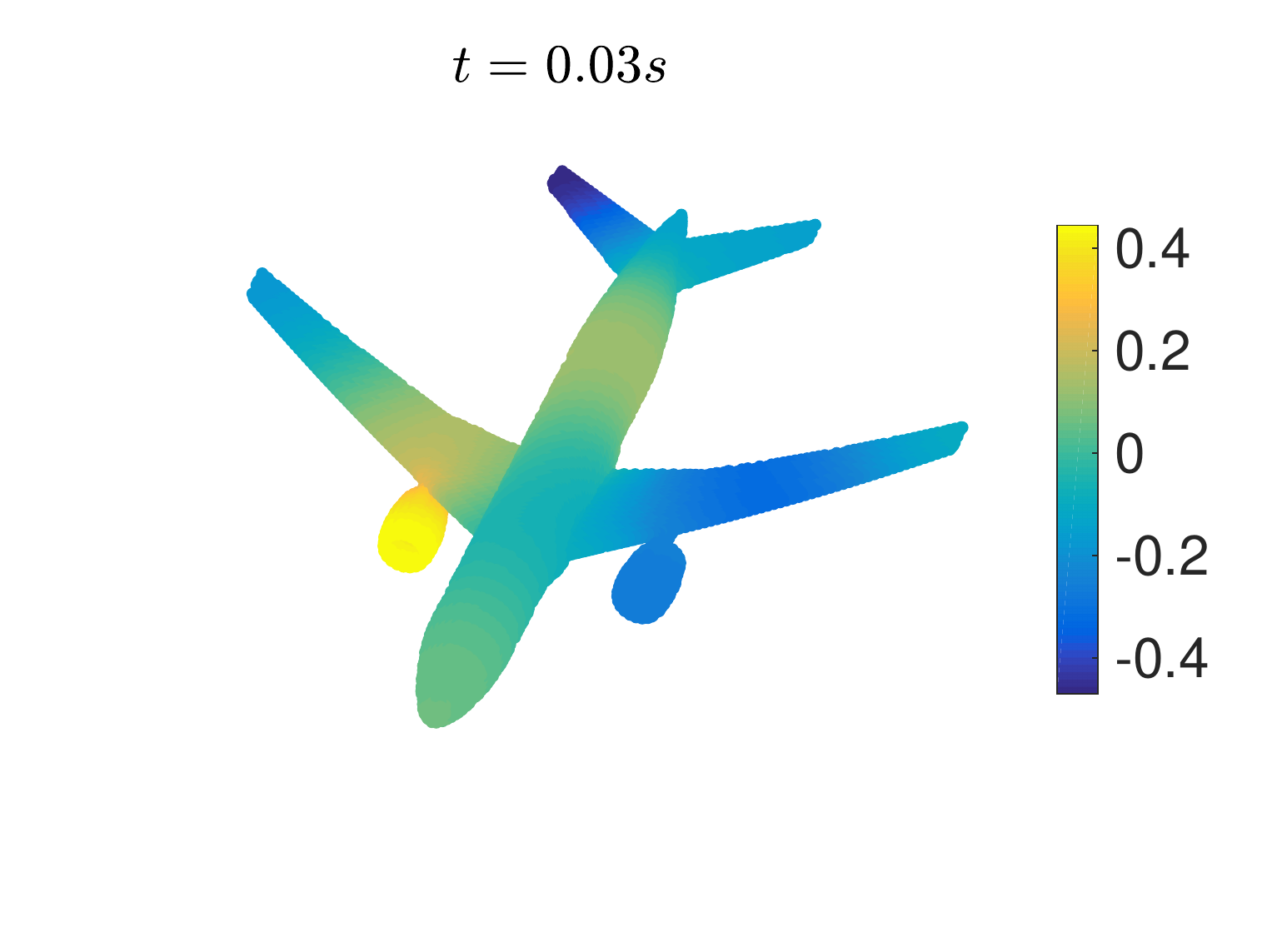}\\
  \includegraphics[width=0.48\textwidth]{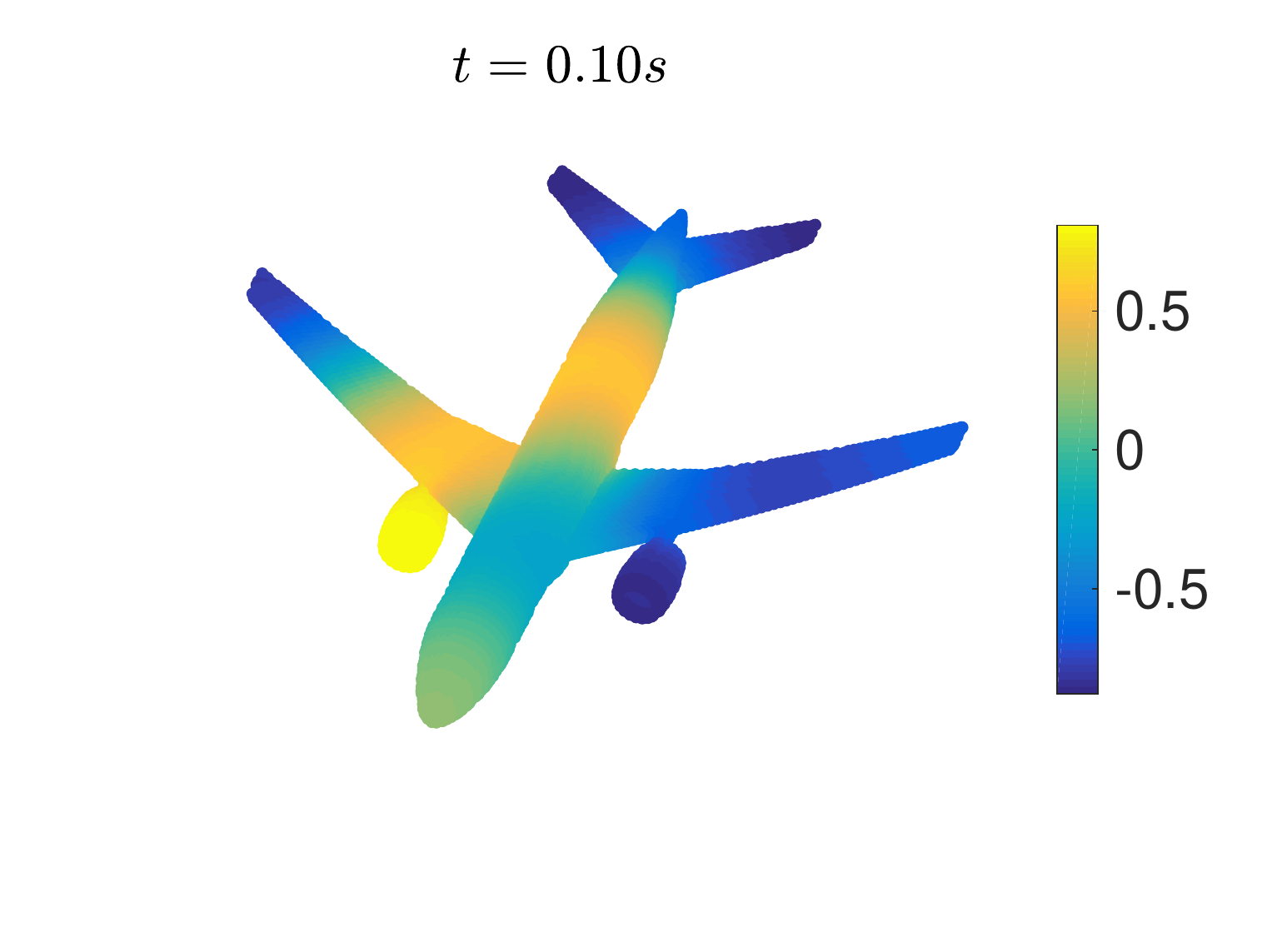}
  \includegraphics[width=0.48\textwidth]{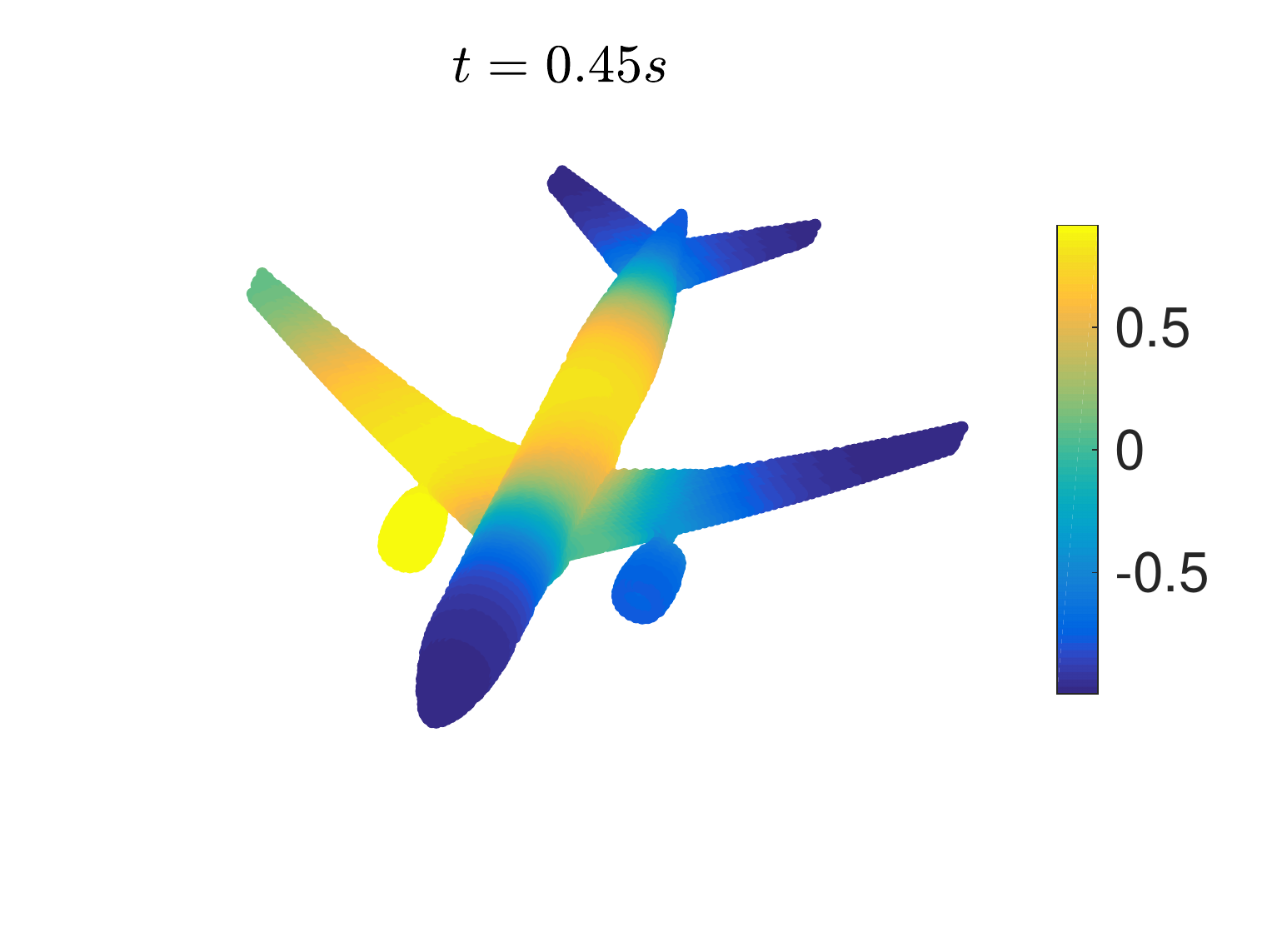}\\
  \includegraphics[width=0.48\textwidth]{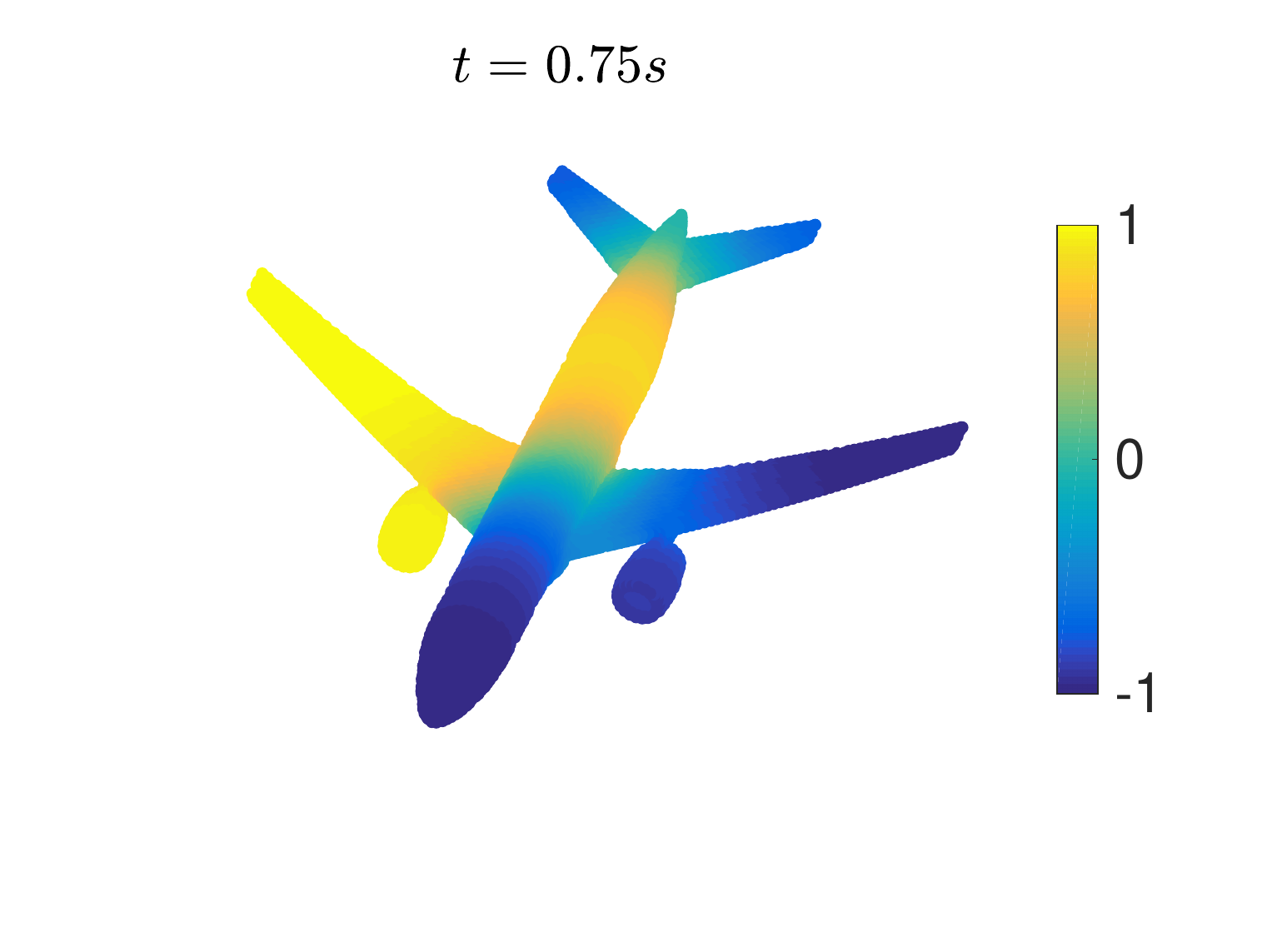}
  \includegraphics[width=0.48\textwidth]{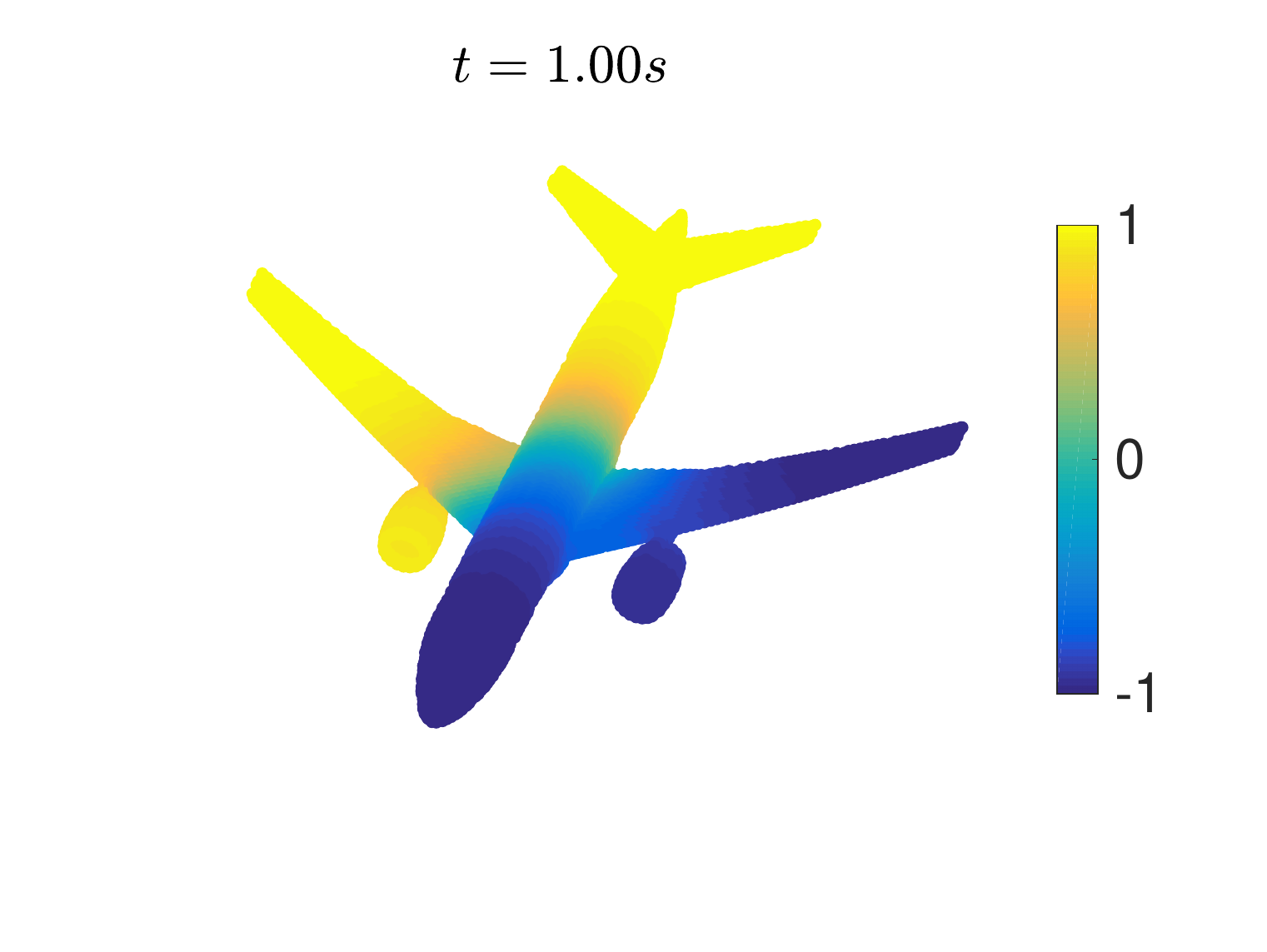}\\    
  \caption{Cahn-Hilliard Equation on the surface of an airplane: evolution of the order parameter $f$ at different times.}%
  % Figures created in CahnHilliard_on_surfaces_plotting.m
  \label{Fig:CahnHilliard_Results}
\end{figure}

This example also illustrates that the present method can handle complex geometries and non-smooth surfaces. The connection between the wings and the fuselage, for instance, has sharp changes. A first order time-integration scheme is used
\begin{align}
	\frac{f^{(n+1)} - f^{(n)}}{\Delta t} &= \frac{1}{Pe} \Delta_{M} \mu^{(n+1)}\,,\\
	\mu^{(n+1)} &= \frac{dg}{df}(f^{(n)}) - Cn^2 \Delta_{M}f^{(n+1)}\,,
\end{align}
which is solved in one large coupled linear system. The evolution of $f$ is shown in Figure~\ref{Fig:CahnHilliard_Results} for $Pe = 1$, $Cn = 0.5$, $\Delta t = 10^{-4}$and $N= 15\,713$ points in the domain. Our simulation results show similar trends to those exhibited in mesh-based solutions \cite{Chen2014,Du2011}. The concentration $f$ first evolves to the range $[-1,1]$, and the separation continues until a steady state is reached.
\section{Conclusion}
\label{sec:Conclusion}

We presented a new meshfree approach to solving PDEs defined on manifolds embedded in $\mathbb{R}^n$. This approach is applicable for manifolds both with and without boundaries. Differential operators are computed directly on the tangent space, and they can be applied directly to function values on the manifold. The discretized domain consists of points only on the manifold, with no bulk discretization around the surface. Consequently, the method scales with the actual dimension of the manifold, and \textit{not} the dimension of the embedding space.

One of the biggest advantage of this method lies in the fact that the differential operators are computed on the tangent space in the same way as regular, volume-based, meshfree methods for the dimension of the tangent space. As a result, all developements in classical volumetric meshfree GFDMs can be directly carried over to manifolds. This was illustrated with different examples, for the optimization of the Laplacian stencil, for treating diffusion with large jumps in the diffusion coefficients, and for discretizing advection terms.

The applicability of this method was shown with numerical examples for a parabolic PDE, an elliptic PDE, and a hyperbolic PDE. Numerical results showed good compatibility with analytical solutions. The numerical simulations also showed that this method can handle unevenly distributed particles, and sharp edges in the point cloud. The simulations showed promising results, and suggest that this method could be extended to solve a wider range of problems on surfaces.

% ==========================================

%\section{Add Somewhere}

%Such projection operators have also been used in the context of surface gradients computed with Radial Basis Functions\cite{Fuselier2013}, however the methods for the Laplacian and Diffusion operators differ significantly. In RBF based projection methods, the Laplacian is computed based on the gradients. Such a computation is infeasible in the context of meshfree GFDMs, as they result in poorly conditioned linear systems \cite{Suchde2018_Thesis}.

%Global coordinates, not local.

%Similar to RBF based methods, this "does not rely on any surface-based metrics or intrinsic
%coordinate  systems,  thus  avoiding  any  coordinate  distortions  or  singularities". Quote from \cite{Fuselier2013}

% =======================================

\section*{Ackowledgements}

The first author would like to thank Dr. Tobias Seifarth for his input on advection methods in meshfree GFDMs, and for sharing his code on the same.

% =======================================

\appendix
\section{Optimized Surface Laplacian Stencils}
\label{App:Laplace}

The Laplace stencil is split into two parts. One which satisfies the consistency conditions of Eq.\,\eqref{Eq:TPL_Consistency}, and one which pushes the stencil towards positivity.
\begin{equation}
\label{Eq:LaplaceSplit}
	c_{ij_{T}}^{\Delta_T} = \sigma_{ij_{T}}^{\Delta_{T}} + \alpha^{\Delta_{T}}d_{ij_{T}}^{\Delta_{T}}\,,
\end{equation}
where the $\sigma_{ij_{T}}^{\Delta_{T}}$ satisfy the consistency conditions, and $d_{ij_{T}}^{\Delta_{T}}$ are used to improve the conditioning of the Laplace stencil. They are given by
\begin{align}
	\sum_{j\in S_i}\sigma_{ij_{T}}^{\Delta_{T}}m_{j_{T}} &= \Delta_{T} m(\vec{x}_i) \qquad \forall m\in\mathcal{P}_{T}\,,\\
	\sigma_{ii}^{\Delta_{T}} &= \tilde{A}_c \,,\\
	\text{min } J_i &= \sum_{j\in S_i} \left( \frac{\sigma_{ij_{T}}^{\Delta_{T}}}{W_{ij_{T}}} \right)^2\,,
\end{align}
where $\tilde{A}_c \notin \{0,1\}$ is some fixed central stencil value. And
\begin{align}
	\sum_{j\in S_i}d_{ij_{T}}^{\Delta_{T}} m_{j_{T}} &= 0 \qquad \forall m\in\mathcal{P}_{T}\,,\label{Eq:DD_dStencil_1}\\
	d_{ii}^{\Delta_{T}} &= 1 \,,\label{Eq:DD_dStencil_2}\\
	\text{min } J_i &= \sum_{j\in S_i} \left( \frac{d_{ij_{T}}^{\Delta_{T}}}{W_{ij_{T}}}  \right)^2\,.\label{Eq:DD_dStencil_3} 
\end{align}
We note that the stencil coefficients $\sigma$ and $d$ can be computed by one minimization with different right hand side vectors \cite[Section A.4]{Suchde2018_Thesis}, and thus, performing this procedure does not increase the computation time significantly.  Now, $\alpha^{\Delta_{T}}$ is computed by minimizing the functional
\begin{equation}
	g^{\Delta_{T}} = \frac{\sum_{j \in S_i} \left(c_{ij_{T}}^{\Delta_{T}} \right)^2}{\left (c_{ii_{T}}^{\Delta_{T}} \right)^2} \,.
\end{equation}
Setting $\frac{\partial g^{\Delta_{T}}}{\partial c_{ii}^{\Delta_{T}}} = 0$, we get
\begin{equation}
	\label{Eq:DD_alpha}
	\alpha^{\Delta_{T}} = \frac{\langle\vec{\sigma},\vec{d} \;\rangle \tilde{A}_c - \langle\vec{\sigma},\vec{\sigma}\;\rangle }{ \langle\vec{\sigma},\vec{d}\;\rangle - \langle\vec{d},\vec{d}\;\rangle \tilde{A}_c}\,,
\end{equation}
where 
\begin{equation}
	\langle \vec{\sigma},\vec{d} \;\rangle = \sum_{j \in S_i} \sigma_{ij_{T}}^{\Delta_{T}} d_{ij_{T}}^{\Delta_{T}}\,.
\end{equation} 
The proof of the above has been shown in \cite{Suchde2018_Thesis} (Section 2.5.5 and Appendix A.4). Using Eq.\,\eqref{Eq:DD_alpha}, we can compute the numerical tangent plane Laplacian according to Eq.\,\eqref{Eq:LaplaceSplit}. This, in turn, gives the numerical surface Laplacian.

\section{Extra Consistency Conditions on the Anisotropic Laplacian Operator}
\label{App:Diffusion}

Here, we prove that the extra consistency conditions Eq.\,\eqref{Eq:DiffExtra1} -- Eq.\,\eqref{Eq:DiffExtra3} imposed on the diffusion operator to deal with large jumps in the diffusion coefficient are valid. To simplify notation, we show the same for the volumetric $2$ dimensional case. And, without loss of generality,  we assume that $\vec{s}$ is along the $x$ direction. They are derived by some algebraic manipulations to improve numerical conditioning. Further, these procedures also avoid the numerical computation of derivatives of $\kappa$, and replace it by the computation of derivatives of $\log(\kappa)$.

Firstly, Eq.\,\eqref{Eq:DiffExtra1} is obtained by the addition of the test function $\frac{1}{\kappa}$, the action of the diffusion operator on which can be given by
\begin{align}
\nabla \cdot \left( \kappa\nabla \left( \frac{1}{\kappa} \right) \right) &= \nabla\cdot \left( - \frac{\nabla \kappa}{\kappa} \right) \,,\\
&= -\nabla \cdot \left( \nabla \left(  \log\kappa     \right)\right)\,, \\
&= -\Delta \left( \log\kappa \right)\,.
\end{align}
Eq.\,\eqref{Eq:DiffExtra2} is obtained by the addition of the test function $\frac{\delta x}{\kappa}$, the action of the diffusion operator on which can be given by
\begin{align}
\nabla \cdot \left( \kappa\nabla \left( \frac{\delta x}{\kappa} \right) \right) &= \nabla\cdot \left( \nabla\left( \delta x \right) - \delta x\frac{\nabla \kappa}{\kappa} \right) \,,\\
&= \Delta \left( \delta x \right) - \nabla\left( \delta x \right)\cdot \nabla \left(\log\kappa \right) - \delta x \Delta \left( \log\kappa \right)\,\\
&= 0 - \frac{\partial}{\partial x} \left( \log\kappa \right) - 0\,,
\end{align}
noting that $\delta x$ evaluated at the central point $i$ is $0$. Lastly, Eq.\,\eqref{Eq:DiffExtra3} is obtained by the addition of the test function $\frac{\delta x ^2}{\kappa}$, the action of the diffusion operator on which can be given by
\begin{align}
\nabla \cdot \left( \kappa\nabla \left( \frac{\delta x^2}{\kappa} \right) \right) &= \nabla\cdot \left( \nabla\left( \delta x \right)^2 - \delta x^2\frac{\nabla \kappa}{\kappa} \right) \,,\\
&= \Delta \left( \delta x \right)^2 - \nabla\left( \delta x \right)^2\cdot \nabla \left(\log\kappa \right) - \delta x^2 \Delta \left( \log\kappa \right)\,\\
&= 2 -0 - 0\,.
\end{align}


\begin{thebibliography}{10}

\bibitem{Alexander1977}
R.~Alexander.
\newblock Diagonally implicit runge–kutta methods for stiff o.d.e.’s.
\newblock {\em SIAM Journal on Numerical Analysis}, 14(6):1006--1021, 1977.

\bibitem{Bertalmio2001}
M.~Bertalm\'{i}o, L.-T. Cheng, S.~Osher, and G.~Sapiro.
\newblock Variational problems and partial differential equations on implicit
  surfaces.
\newblock {\em Journal of Computational Physics}, 174(2):759 -- 780, 2001.

\bibitem{Cahn1958}
J.~W. Cahn and J.~E. Hilliard.
\newblock Free energy of a nonuniform system. i. interfacial free energy.
\newblock {\em The Journal of Chemical Physics}, 28(2):258--267, 1958.

\bibitem{Chen2014}
S.~Chen and J.~Wu.
\newblock Discrete conservation laws on curved surfaces {II}: A dual approach.
\newblock {\em SIAM Journal on Scientific Computing}, 36(4):A1813--A1830, 2014.

\bibitem{ChenThesis_CPM}
Y.~Chen.
\newblock {\em Geometric Multigrid and Closest Point Methods for Surfaces and
  General Domains}.
\newblock PhD thesis, St Anne's College, University of Oxford, Oxford, 2015.

\bibitem{Chu2018}
J.~Chu and R.~Tsai.
\newblock Volumetric variational principles for a class of partial differential
  equations defined on surfaces and curves.
\newblock {\em Research in the Mathematical Sciences}, 5(2):19, 2018.

\bibitem{Davydov2011}
O.~Davydov and D.~T. Oanh.
\newblock On the optimal shape parameter for gaussian radial basis function
  finite difference approximation of the poisson equation.
\newblock {\em Computers \& Mathematics with Applications}, 62(5):2143 -- 2161,
  2011.

\bibitem{Demanet2006}
L.~Demanet.
\newblock Painless, highly accurate discretizations of the laplacian on a
  smooth manifold.
\newblock Technical report, Stanford University, 2006.

\bibitem{Diewald2000}
U.~Diewald, T.~Preusser, and M.~Rumpf.
\newblock Anisotropic diffusion in vector field visualization on euclidean
  domains and surfaces.
\newblock {\em IEEE Transactions on Visualization and Computer Graphics},
  6(2):139--149, Apr 2000.

\bibitem{Drumm2008}
C.~Drumm, S.~Tiwari, J.~Kuhnert, and H.-J. Bart.
\newblock Finite pointset method for simulation of the liquid - liquid flow
  field in an extractor.
\newblock {\em Computers \& Chemical Engineering}, 32(12):2946 -- 2957, 2008.

\bibitem{Du2003}
Q.~Du, M.~D. Gunzburger, and L.~Ju.
\newblock Voronoi-based finite volume methods, optimal voronoi meshes, and pdes
  on the sphere.
\newblock {\em Computer Methods in Applied Mechanics and Engineering},
  192(35):3933 -- 3957, 2003.
  
\bibitem{Du2006}
Q.~Du and L.~Ju.
\newblock Finite volume methods on spheres and spherical centroidal voronoi
  meshes.
\newblock {\em SIAM Journal on Numerical Analysis}, 43(4):1673--1692, 2005.

\bibitem{Du2011}
Q.~Du, L.~Ju, and L.~Tian.
\newblock Finite element approximation of the cahn-hilliard equation on
  surfaces.
\newblock {\em Computer Methods in Applied Mechanics and Engineering},
  200(29):2458 -- 2470, 2011.

\bibitem{Dziuk2013}
G.~Dziuk and C.~M. Elliott.
\newblock Finite element methods for surface pdes.
\newblock {\em Acta Numerica}, 22:289--396, 2013.

\bibitem{Ellsiepen1999}
P.~Ellsiepen.
\newblock {\em Zeit-und ortsadaptive Verfahren angewandt auf Mehrphasenprobleme
  por{\"o}ser Medien}.
\newblock PhD thesis, Stuttgart University, Stuttgart, 1999.

\bibitem{Fan2018}
C.-M. Fan, C.-N. Chu, B.~Šarler, and T.-H. Li.
\newblock Numerical solutions of waves-current interactions by generalized
  finite difference method.
\newblock {\em Engineering Analysis with Boundary Elements}, 2018.

\bibitem{Floater2005}
M.~S. Floater and K.~Hormann.
\newblock Surface parameterization: a tutorial and survey.
\newblock In N.~A. Dodgson, M.~S. Floater, and M.~A. Sabin, editors, {\em
  Advances in Multiresolution for Geometric Modelling}, pages 157--186, Berlin,
  Heidelberg, 2005. Springer Berlin Heidelberg.

\bibitem{Flyer2014}
N.~Flyer, G.~B. Wright, and B.~Fornberg.
\newblock Radial basis function-generated finite differences: A mesh-free
  method for computational geosciences.
\newblock In W.~Freeden, M.~Z. Nashed, and T.~Sonar, editors, {\em Handbook of
  Geomathematics}, pages 1--30, Berlin, Heidelberg, 2014. Springer Berlin
  Heidelberg.

\bibitem{Froese2018}
B.~D. Froese.
\newblock Meshfree finite difference approximations for functions of the
  eigenvalues of the hessian.
\newblock {\em Numerische Mathematik}, 138(1):75--99, Jan 2018.

\bibitem{Fuselier2013}
E.~J. Fuselier and G.~B. Wright.
\newblock A high-order kernel method for diffusion and reaction-diffusion
  equations on surfaces.
\newblock {\em Journal of Scientific Computing}, 56(3):535--565, Sep 2013.

\bibitem{Gavete2017}
L.~Gavete, F.~Ure{\~n}a, J.~Benito, A.~Garc\'{i}a, M.~Ure{\~n}a, and E.~Salete.
\newblock Solving second order non-linear elliptic partial differential
  equations using generalized finite difference method.
\newblock {\em Journal of Computational and Applied Mathematics}, 318:378 --
  387, 2017.
\newblock Computational and Mathematical Methods in Science and Engineering
  CMMSE-2015.

\bibitem{Gera2017}
P.~Gera and D.~Salac.
\newblock Cahn-hilliard on surfaces: A numerical study.
\newblock {\em Applied Mathematics Letters}, 73:56 -- 61, 2017.

\bibitem{Jefferies2015}
A.~Jefferies, J.~Kuhnert, L.~Aschenbrenner, and U.~Giffhorn.
\newblock Finite pointset method for the simulation of a vehicle travelling
  through a body of water.
\newblock In M.~Griebel and A.~M. Schweitzer, editors, {\em Meshfree Methods
  for Partial Differential Equations VII}, pages 205--221, Cham, 2015. Springer
  International Publishing.

\bibitem{Katz2010}
A.~Katz and A.~Jameson.
\newblock Meshless scheme based on alignment constraints.
\newblock {\em AIAA journal}, 48(11):2501--2511, 2010.

\bibitem{Kim2016}
J.~Kim, S.~Lee, Y.~Choi, S.-M. Lee, and D.~Jeong.
\newblock Basic principles and practical applications of the cahn--hilliard
  equation.
\newblock {\em Mathematical Problems in Engineering}, 2016, 2016.

\bibitem{Chiu2012}
E.~Kwan-yu Chiu, Q.~Wang, R.~Hu, and A.~Jameson.
\newblock A conservative mesh-free scheme and generalized framework for
  conservation laws.
\newblock {\em SIAM Journal on Scientific Computing}, 34(6):A2896--A2916, 2012.

\bibitem{Lai2013}
R.~Lai, J.~Liang, and H.~Zhao.
\newblock A local mesh method for solving pdes on point clouds.
\newblock {\em Inverse Problems \& Imaging}, 7(3), 2013.

\bibitem{Liang2012}
J.~Liang, R.~Lai, T.~W. Wong, and H.~Zhao.
\newblock Geometric understanding of point clouds using laplace-beltrami
  operator.
\newblock In {\em Computer Vision and Pattern Recognition (CVPR), 2012 IEEE
  Conference on}, pages 214--221. IEEE, 2012.

\bibitem{Liang2013}
J.~Liang and H.~Zhao.
\newblock Solving partial differential equations on point clouds.
\newblock {\em SIAM Journal on Scientific Computing}, 35(3):A1461--A1486, 2013.

\bibitem{Luo2016}
M.~{Luo}, C.~G. {Koh}, W.~{Bai}, and M.~{Gao}.
\newblock {A particle method for two-phase flows with compressible air pocket}.
\newblock {\em International Journal for Numerical Methods in Engineering},
  108:695--721, Nov. 2016.

\bibitem{Marz2012}
T.~M\"arz and C.~B. Macdonald.
\newblock Calculus on surfaces with general closest point functions.
\newblock {\em SIAM Journal on Numerical Analysis}, 50(6):3303--3328, 2012.

\bibitem{Milewski2012}
S.~Milewski.
\newblock Meshless finite difference method with higher order
  approximation---applications in mechanics.
\newblock {\em Archives of Computational Methods in Engineering}, 19(1):1--49,
  2012.

\bibitem{Mitra2003}
N.~J. Mitra and A.~Nguyen.
\newblock Estimating surface normals in noisy point cloud data.
\newblock In {\em Proceedings of the Nineteenth Annual Symposium on
  Computational Geometry}, SCG '03, pages 322--328, New York, NY, USA, 2003.
  ACM.

\bibitem{Moller2007}
A.~M{\"o}ller and J.~Kuhnert.
\newblock Simulation of the glass flow inside a floating process / {S}imulation
  de l'{\'e}coulement du verre dans le proc{\'e}d{\'e} float.
\newblock {\em Revue Verre}, 13(5):28--30, 2007.

\bibitem{Mongillo2011}
M.~Mongillo.
\newblock Choosing basis functions and shape parameters for radial basis
  function methods.
\newblock {\em SIAM Undergraduate Research Online}, 4:190--209, 2011.

\bibitem{Myers2002}
T.~G. {Myers}, J.~P.~F. {Charpin}, and S.~J. {Chapman}.
\newblock {The flow and solidification of a thin fluid film on an arbitrary
  three-dimensional surface}.
\newblock {\em Physics of Fluids}, 14:2788--2803, Aug. 2002.

\bibitem{Novak2007}
I.~L. Novak, F.~Gao, Y.-S. Choi, D.~Resasco, J.~C. Schaff, and B.~M.
  Slepchenko.
\newblock Diffusion on a curved surface coupled to diffusion in the volume:
  Application to cell biology.
\newblock {\em Journal of Computational Physics}, 226(2):1271 -- 1290, 2007.

\bibitem{Olshanskii2009}
M.~A. Olshanskii, A.~Reusken, and J.~Grande.
\newblock A finite element method for elliptic equations on surfaces.
\newblock {\em SIAM Journal on Numerical Analysis}, 47(5):3339--3358, 2009.

\bibitem{Persson2004}
P.-O. Persson and G.~Strang.
\newblock A simple mesh generator in matlab.
\newblock {\em SIAM Review}, 46:2004, 2004.

\bibitem{Petronetto2013}
F.~Petronetto, A.~Paiva, E.~S. Helou, D.~E. Stewart, and L.~G. Nonato.
\newblock Mesh-free discrete laplace-beltrami operator.
\newblock {\em Comput. Graph. Forum}, 32(6):214--226, Sept. 2013.

\bibitem{Piret2012}
C.~Piret.
\newblock The orthogonal gradients method: A radial basis functions method for
  solving partial differential equations on arbitrary surfaces.
\newblock {\em Journal of Computational Physics}, 231(14):4662 -- 4675, 2012.

\bibitem{Praveen2007}
C.~Praveen and S.~M. Deshpande.
\newblock Kinetic meshless method for compressible flows.
\newblock {\em International Journal for Numerical Methods in Fluids},
  55(11):1059--1089, 2007.

\bibitem{Ratz2006}
A.~R{\"a}tz and A.~Voigt.
\newblock Pde's on surfaces---a diffuse interface approach.
\newblock {\em Commun. Math. Sci.}, 4(3):575--590, 09 2006.


\bibitem{Edgar2017}
E.~O. Res{\'e}ndiz-Flores, J.~Kuhnert, and F.~R. Saucedo-Zendejo.
\newblock Application of a generalized finite difference method to mould
  filling process.
\newblock {\em European Journal of Applied Mathematics}, page 1–20, 2017.

\bibitem{Ruuth2008}
S.~J. Ruuth and B.~Merriman.
\newblock A simple embedding method for solving partial differential equations
  on surfaces.
\newblock {\em Journal of Computational Physics}, 227(3):1943 -- 1961, 2008.

\bibitem{SeiboldThesis}
B.~Seibold.
\newblock {\em M-Matrices in Meshless Finite Difference Methods}.
\newblock PhD thesis, Kaiserslautern University, 2006.

\bibitem{Seibold2008}
B.~Seibold.
\newblock Minimal positive stencils in meshfree finite difference methods for
  the poisson equation.
\newblock {\em Computer Methods in Applied Mechanics and Engineering},
  198(3-4):592 -- 601, 2008.

\bibitem{SeifarthThesis}
T.~Seifarth.
\newblock {\em Numerische Algortihmen f{\"u}r gitterfreie Methoden zur
  L{\"o}sung von Transportproblemen}.
\newblock PhD thesis, University of Kassel, Kassel, 2017.

\bibitem{Simonenko2014}
S.~Simonenko, V.~Bayona, and M.~Kindelan.
\newblock Optimal shape parameter for the solution of elastostatic problems
  with the rbf method.
\newblock {\em Journal of Engineering Mathematics}, 85(1):115--129, Apr 2014.

\bibitem{Suchde2018_Thesis}
P.~Suchde.
\newblock {\em Conservation and Accuracy in Meshfree Generalized Finite
  Difference Methods}.
\newblock PhD thesis, University of Kaiserslautern, Kaiserslautern, Germany,
  2018.

\bibitem{Suchde2017_CCC}
P.~Suchde, J.~Kuhnert, S.~Schr\"oder, and A.~Klar.
\newblock A flux conserving meshfree method for conservation laws.
\newblock {\em International Journal for Numerical Methods in Engineering},
  112(3):238--256, 2017.

\bibitem{Suchde2018_INSE}
P.~Suchde, J.~Kuhnert, and S.~Tiwari.
\newblock On meshfree {GFDM} solvers for the incompressible {N}avier--{S}tokes
  equations.
\newblock {\em Computers \& Fluids}, 165:1 -- 12, 2018.

\bibitem{Tramecon2013}
A.~Tramecon and J.~Kuhnert.
\newblock Simulation of advanced folded airbags with {VPS-PAMCRASH/FPM}:
  {D}evelopment and validation of turbulent flow numerical simulation
  techniques applied to curtain bag deployments.
\newblock In {\em SAE Technical Paper}, Warrendale, PA, USA, 2013. SAE
  International.

\bibitem{Trask2017}
N.~Trask, M.~Perego, and P.~B. Bochev.
\newblock A high-order staggered meshless method for elliptic problems.
\newblock {\em SIAM J. Scientific Computing}, 39(2):A479--A502, 2017.

\bibitem{Turk1991}
G.~Turk.
\newblock Generating textures on arbitrary surfaces using reaction-diffusion.
\newblock {\em SIGGRAPH Comput. Graph.}, 25(4):289--298, July 1991.

\bibitem{BiCGSTAB}
H.~A. van~der Vorst.
\newblock Bi-cgstab: A fast and smoothly converging variant of bi-cg for the
  solution of nonsymmetric linear systems.
\newblock {\em SIAM Journal on Scientific and Statistical Computing},
  13(2):631--644, 1992.

\bibitem{Vassberg2008}
J.~Vassberg, M.~Dehaan, M.~Rivers, and R.~Wahls.
\newblock Development of a common research model for applied {CFD} validation
  studies.
\newblock In {\em 26th AIAA Applied Aerodynamics Conference}, page 6919, 2008.

\bibitem{IngridThesis_CPM}
I.~von Glehn.
\newblock {\em A closest point penalty method for evolution equations on
  surfaces}.
\newblock PhD thesis, Oriel College, University of Oxford, Oxford, 2014.

\bibitem{Yoon2014}
Y.-C. Yoon and J.-H. Song.
\newblock Extended particle difference method for weak and strong discontinuity
  problems: part i. derivation of the extended particle derivative
  approximation for the representation of weak and strong discontinuities.
\newblock {\em Computational Mechanics}, 53(6):1087--1103, Jun 2014.

\end{thebibliography}
\end{document}